\documentclass[preprint]{elsarticle}
\usepackage[T1]{fontenc}
\usepackage[utf8]{inputenc}
\usepackage[letterpaper, margin=1in]{geometry}
\usepackage[colorlinks]{hyperref}
\usepackage{subfigure}
\usepackage{array}
\usepackage{amsmath,amssymb,amsfonts,xcolor,amsthm,dsfont,,mathabx}
\usepackage{mathtools, physics}
\usepackage{placeins}
\usepackage{csquotes}
\usepackage[capitalise, nameinlink, noabbrev, nosort]{cleveref}
\usepackage{doi}

\usepackage{latexsym}
\usepackage{euscript}
\usepackage{tocvsec2}
\usepackage{etoolbox,caption}
\usepackage{float}
\usepackage{marginnote}
\usepackage{epsfig,tikz}
\usepackage{tikz-cd}%commutative diagrams
\usetikzlibrary{patterns}
\usetikzlibrary{arrows.meta}
\usetikzlibrary{matrix}

\usepackage{enumitem}
\makeatletter
\newcommand{\mylabel}[2]{#2\def\@currentlabel{#2}\label{#1}}
\makeatother

\newtheorem*{theorem-no}{Theorem}

\numberwithin{pulse}{section}  % all thms, etc set from pulse

\makeatletter\g@addto@macro\bfseries{\boldmath}\makeatother

\makeatletter
\let\origsection\section
\renewcommand\section{\@ifstar{\starsection}{\nostarsection}}
\newcommand\sectionspace{\vspace{0.5ex}}
\newcommand\nostarsection[1]{\sectionspace\origsection{#1}\sectionspace}
\newcommand\starsection[1]{\sectionspace\origsection*{#1}\sectionspace}
\makeatother

\setlist[enumerate]{font=\normalfont}
\crefname{enumi}{}{}

\crefname{case}{case}{cases}
\creflabelformat{case}{#2#1#3}

\usepackage[hang, flushmargin]{footmisc}
\usepackage[english]{isodate}
\usepackage{mathrsfs}

\usepackage[capitalise, nameinlink, noabbrev, nosort]{cleveref}

\allowdisplaybreaks

\setenumerate{listparindent=\parindent}

\crefname{page}{page}{pages}

\numberwithin{equation}{section}
\crefname{equation}{equation}{equations}

\crefname{inequality}{inequality}{inequalities}
\creflabelformat{inequality}{#2(#1)#3}

\crefname{condition}{condition}{conditions}
\creflabelformat{condition}{#2(#1)#3}

\newtheorem{theorem}{Theorem}[section]

\newtheorem{thm}[theorem]{Theorem}
\crefname{thm}{Theorem}{Theorems}

\newtheorem{lem}[theorem]{Lemma}
\crefname{lem}{Lemma}{Lemmas}

\newtheorem{prop}[theorem]{Proposition}
\crefname{prop}{Proposition}{Propositions}

\newtheorem{cor}[theorem]{Corollary}
\crefname{cor}{Corollary}{Corollaries}

\crefname{assumption}{Assumption}{Assumptions}

\theoremstyle{definition}

\newtheorem{defn}[theorem]{Definition}
\crefname{defn}{Definition}{Definitions}

\newtheorem{rmk}[theorem]{Remark}
\crefname{rmk}{Remark}{Remarks}

\newtheorem{example}[theorem]{Example}
\crefname{example}{Example}{Examples}

\crefname{question}{Question}{Questions}

\newtheorem{claim}[theorem]{Claim}
\newtheorem*{claim*}{Claim}%creates theorem, proposition, conjecture, corollary, claim and lemma environments

\newcommand{\set}[1]{\left\{ #1 \right\}} %Set brackets
\newcommand{\innprod}[1]{\left \langle #1 \right \rangle} %inner product braces
\newcommand{\nnorm}[1]{\|#1\|} %norm brackets
\newcommand{\of}[1]{\left( #1 \right)} %parenthesis
\newcommand{\supp}{\operatorname{supp}}
\newcommand{\Hil}{\mathcal{H}}
\newcommand{\fH}{\mathfrak{H}}

\newcommand{\F}{\mathbb{F}}

\newcommand{\G}{\mathcal{G}}
\newcommand{\N}{\mathbb{N}}
\newcommand{\C}{\mathbb{C}}

\newcommand{\unsp}{\G^{(0)}}
\newcommand{\prodsp}{\G^{(2)}} 
\newcommand{\ray}{\text{ray}}
\newcommand{\bxi}{\boldsymbol{\xi}}
\newcommand{\boldeta}{\boldsymbol{\eta}}

\newcommand{\infN}{\N^*}
\newcommand{\calpi}{\varpi}

\renewcommand\paragraph[1]{\par\vspace{1em}\noindent\textbf{#1}}

\makeatletter
\def\ps@pprintTitle{%
  \let\@oddhead\@empty
  \let\@evenhead\@empty
  \let\@oddfoot\@empty
  \let\@evenfoot\@oddfoot
}
\makeatother

\begin{document}
\begin{frontmatter}
\title{Fourier Algebras of HLS Groupoids: Lifts, Amenability, and Spectrum}

\author[mysecondaryaddress]{Joseph DeGaetani}
\author[mymainaddress]{Mahya Ghandehari\corref{mycorrespondingauthor}}
\cortext[mycorrespondingauthor]{Corresponding author}
\ead{mahya@udel.edu}
\address[mysecondaryaddress]{Department of Mathematics, State University of New York at Oswego, 
Oswego, NY, 13126, USA}
\address[mymainaddress]{Department of Mathematical Sciences, University of Delaware, Newark, DE, 19716, USA}
%
%\date{16 Feb 2023}
%
\begin{abstract}
%In \cite{PatersonPaper}, Paterson introduced a definition of the Fourier algebra of a locally compact groupoid based on continuous Hilbert bundles and established its basic properties. 
%
We investigate the Fourier algebra of a certain class of groupoids, called the HLS groupoids, within Paterson's continuous framework for locally compact groupoids. We develop the basic harmonic analysis of HLS groupoids by obtaining concrete descriptions of their topology, representations, and Fourier algebra. 
We then study the restriction maps from the Fourier algebra $A(\G)$ of an HLS groupoid $\G$ to the Fourier algebras of its isotropy subgroups, and construct explicit norm-preserving lifts for these restriction maps.
As applications, we describe the norm of every element of $A(\mathcal{G})$ in terms of the norms of its restrictions to the isotropy subgroups.  We prove an analogue of Leptin's theorem characterizing when $A(\mathcal{G})$ has a bounded approximate identity. We also show that the Gelfand spectrum of $A(\mathcal{G})$ is homeomorphic to $\mathcal{G}$. 
Finally, using abstract arguments that, unlike our methods for HLS groupoids, do not yield explicit constructions, we show that norm-preserving lifts exist and that the Gelfand spectrum of the Fourier algebra is homeomorphic to the groupoid for every \'{e}tale groupoid. This generalizes Eymard's theorem and answers a question of Paterson for \'{e}tale groupoids.
\end{abstract}
%*46J99 - Commutative Banach algebras and commutative topological algebras
%*43A30 - Fourier and Fourier--Stieltjes transforms on nonabelian groups and on semigroups, etc.
%*43A32 - Other transforms and operators of Fourier type
%*46J10 - Banach algebras of continuous functions, function algebras
%*22A22 - Topological groupoids (including differentiable and Lie groupoids) 
\end{frontmatter}

%%%%%%%%%%%%%%%%%%%%%%%%%%%%%%%%%%%%%%%%%%%%%%%%%%%%%%%%%%%%%%%%%%%%%%%%%%%%%%%%%%%%%%%%%%%%%%%%
\section{Introduction}
In \cite{Eymard}, Eymard introduced the Fourier and Fourier--Stieltjes algebras of a locally compact group $G$ and initiated their systematic study. Denoted by $A(G)$ and $B(G)$, respectively, these Banach algebras play a fundamental role in abstract harmonic analysis. When $G$ is Abelian, the Fourier transform identifies $A(G)$ and $B(G)$ with the group algebra and the measure algebra of the dual group $\widehat{G}$, respectively. In the non-Abelian setting, $B(G)$ is naturally identified with the continuous dual of the full group $C^*$-algebra, while $A(G)$ is the predual of the group von Neumann algebra. 
The Fourier and Fourier--Stieltjes algebras encode the properties of the underlying group. 
For instance, Eymard proved that the Gelfand spectrum of $A(G)$, namely the set of nonzero multiplicative linear functionals on $A(G)$, is homeomorphic to $G$ via point evaluation.

Analogous function algebras, denoted $A(\G)$ and $B(\G)$, have been defined when the underlying structure is not a group but a groupoid $\G$. Groupoids, first introduced by Brandt in \cite{Brandt}, are algebraic structures that generalize groups. They require only a partially defined product and do not necessarily have a unique identity element. Instead, they possess a unit space $\unsp$, whose elements serve as local identities.
The theories of the Fourier and Fourier--Stieltjes algebras for locally compact groupoids were developed in recent years from multiple perspectives \cite{DeG, Oty,  Paterson-2003, PatersonPaper, Ramsay-Martin-97, Renault}. These theories largely differ in the choice between requiring continuity or measurability, and the choice of norm. 
In the special case where $\G$ is a locally compact group, all of these theories reduce to the existing definitions for the Fourier and Fourier--Stieltjes algebras of groups. 
In this paper, we study the Fourier algebra of locally compact groupoids, adopting the approach of \cite{PatersonPaper} by developing a continuous, rather than a measurable, theory.
In this approach, a representation of a groupoid is defined to be a family of unitary maps between Hilbert spaces from a continuous Hilbert bundle, with coefficient functions being obtained using continuous sections of the Hilbert field. 

While the Fourier and Fourier--Stieltjes algebras of locally compact groups have been studied extensively, their counterparts for locally compact groupoids remain comparatively unexplored. In particular, relatively little is known about their properties as function algebras or about the extent to which they reflect the structure of the underlying groupoid.
The primary focus of this manuscript is the Fourier algebra of a class of groupoids known as HLS groupoids, continuing the investigation initiated in \cite{DeG}. Motivated by analogous results for locally compact groups (see, for example, \cite{Cowling, MahyaThesisPaper, Herz, HewittII, Leptin}), we studied the restriction maps for the Fourier and Fourier--Stieltjes algebras of locally compact groupoids in \cite{DeG}, establishing several general functorial properties.
In particular, we proved that for every \'{e}tale groupoid the restriction map from the Fourier algebra to the Fourier algebra of any isotropy subgroup is surjective, while the analogous statement for the Fourier--Stieltjes algebra fails in general. Complete descriptions of the Fourier and Fourier--Stieltjes algebras were also obtained for group bundles with discrete unit space.

Harmonic analysis on group bundles has been developed mainly in the Abelian setting through Pontryagin duality, fibrewise Fourier transforms, and dual Haar systems \cite{GoehleGroupBundle,MuhlyRenaultWilliamsIII,RenaultDualHaar,ShirinkalamMedghalchiAmini}, and in the compact setting through Peter--Weyl-type representation theory \cite{EdekoJamneshanKreidler}. By contrast, little is known about the Fourier and Fourier--Stieltjes algebras of non-discrete, non-Abelian group bundles.
HLS groupoids therefore provide a natural setting in which to investigate how the topology of a group bundle influences its Fourier algebra.
Although HLS groupoids are group bundles, their unit space is no longer discrete, and the topology couples the finite fibers with the infinite fiber in a nontrivial way. Consequently, the techniques developed in \cite{DeG} for discrete group bundles no longer apply here, and several basic structural questions concerning the Fourier algebra become substantially more delicate. 
On the other hand, the topology of HLS groupoids is sufficiently explicit to allow a detailed analysis. 
At the same time, HLS groupoids already exhibit several interesting phenomena absent in the settings of groups and discrete group bundles, including the failure of surjectivity of the Fourier--Stieltjes restriction map \cite{DeG} and the failure of the converse to the groupoid amenability theorem \cite{Willett}.

\subsection{Main contributions} 
Our results fall into two categories. Most of the paper concerns HLS groupoids.
The topology of an HLS groupoid admits a particularly clear description that simplifies the associated Fourier theory. As a result, we obtain concrete formulas and a structural description of the algebra, rather than existence statements alone. In the final section, we discuss the Fourier algebra of general \'{e}tale groupoids. The general arguments are abstract and do not recover the finer information available for HLS groupoids.

\medskip
\noindent\textbf{Fourier algebra of HLS groupoids.}
Taking advantage of the special topological structure of HLS groupoids, we obtain new insights into their Fourier algebras and establish useful bounds for the Fourier norm (\Cref{subsec:fourier-hls}).
In particular, we identify the continuous sections of the left regular Hilbert bundle with an explicit class of continuous functions on $\G$ (\Cref{cor-cts-sec}). This concrete picture of the left regular representation underlies many subsequent results.

%lift
In our previous work \cite{DeG}, we proved that the restriction map from the Fourier algebra of an \'{e}tale groupoid to an isotropy subgroup is surjective. 
Here we  strengthen that result for HLS groupoids by constructing explicit norm-preserving lifts of Fourier algebra functions from every fiber (\Cref{thm:HLSSurjective}). 
The construction preserves positive definiteness, and it allows elements on finitely many fibers to be lifted simultaneously.

In \Cref{sec:application}, we discuss two applications of the lifting result (\Cref{thm:HLSSurjective}).
First, the explicit lifts let us describe the norm of an element in $A(\G)$ completely in terms of its fibers
(\Cref{cor:norm-convergence-fibers-kernel,cor:norm-A(G)-limit-sections}). Consequently, $A(\G)$ embeds isometrically into the $\ell^\infty$-direct sum of the Fourier algebras of the isotropy subgroups of $\G$.
Second, we prove a Leptin-type theorem for HLS groupoids: The Fourier algebra $A(\G)$ has a bounded approximate identity if and only if the isotropy subgroup at infinity is amenable, and in that case the approximate identity can be chosen to be bounded by $1$  (\Cref{cor:BAI}). Thus a Banach-algebraic property of $A(\G)$ is governed entirely by the single infinite fiber.

Finally, we investigate the relationship between the Fourier algebra of an HLS groupoid and the underlying groupoid through its Gelfand spectrum. We prove that every nonzero multiplicative linear functional on the Fourier algebra is given by point evaluation, thereby showing that the Gelfand spectrum of the Fourier algebra of an HLS groupoid is naturally homeomorphic to the groupoid itself (\Cref{prop:spectrum}).

\medskip
\noindent\textbf{Beyond HLS groupoids.}
In \Cref{section:beyond-HLS}, we show that two of our results hold for every \'{e}tale groupoid. Norm-preserving lifts exist for the restriction maps onto isotropy subgroups (\Cref{prop:etale-lift}), and the Gelfand spectrum of the Fourier algebra is homeomorphic to the groupoid (\Cref{prop-spec-etale}). The latter extends Eymard's classical theorem for locally compact groups to \'{e}tale groupoids, and answers a question of Paterson \cite[Section~6]{PatersonPaper} for the \'{e}tale setting; the question remains open in general. These arguments, however, are abstract: the lifts are obtained from a non-constructive extension of vectors to continuous sections, and they come with no explicit formula. Our results on the fiberwise structure of the norm, the limiting behavior of the norm, and the bounded approximate identity all rely on the explicit form of the HLS lifts and on the HLS topology. We do not know analogues of these results for general \'{e}tale groupoids. This is the main reason the HLS class is worth studying in its own right.

%%%%%%%%%%%%%%%%%%%%%%%%%%%%%%%%%%%%%%%%%%%%%%%%%%%%%%%%%%%%%%%%%%%%%%%%%%%%%%%%%%%%%%%%%%%%%%%%
%%%%%%%%%%%%%%%%%%%%%%%%%%%%%%%%%%%%%%%%%%%%%%%%%%%%%%%%%%%%%%%%%%%%%%%%%%%%%%%%%%%%%%%%%%%%%%%%
\section{Notations and background}
\label{sec:NotationsAndBackground}

%%%%%%%%%%%%%%%%%%%%%%%%%%%%%%%%%%%%%%%%%%%%%%%%%%%%%%%%%%%%%%%%%%%%%%%%%%%%%%%%%%%%%%%%%%%%%%%%
A \emph{groupoid} is a set $\G$ together with a subset $\prodsp\subseteq \G\times \G$, a product map $\prodsp \to \G$, $(a,b)\mapsto ab$ and an involutive inverse map $\G\to \G$, $a\mapsto a^{-1}$ that satisfies the following:
\begin{enumerate}
    \item if $(a,b),(b,c)\in \prodsp$, then $(ab,c),(a,bc)\in \G^{(2)}$ and $(ab)c=a(bc);$
    \item $(b,b^{-1})\in \prodsp$ for all $b\in \G$, and if $(a,b)\in \prodsp$, then $a^{-1}(ab)=b$, and $(ab)b^{-1}=a.$
\end{enumerate}
Groupoids come equipped with \emph{range} and \emph{source} maps, denoted $r$ and $s$ respectively which are defined as 
$$r(a):=aa^{-1} \; \text{ and }\; s(a):=a^{-1}a.$$
These maps can be used to define the \emph{unit space} as $\unsp:=r(\G)=\set{aa^{-1}:a\in \G}$. Clearly, $\unsp=s(\G)$ as well. 
For $u\in \unsp$, we set $\G^{u}:=r^{-1}\of{\set{u}}$, $\G_u:=s^{-1}\of{\set{u}}$, and $\G_u^v:=\G_u\cap \G^v$. These fibers do not typically form algebraic objects in their own right, however, the set $\G_u^u$ always forms a group, called the \emph{isotropy subgroup at $u$}.

An important class of groupoids is a disjoint union of groups, known as a \emph{group bundle}. 
Namely, for a collection $\{G_i\}_{i\in I}$ of groups, define $\G:=\sqcup_{i\in I} (\{i\}\times G_i)$. 
Elements of $\G$ are composable if they have identical indices in their first coordinates. The product and the inverse are defined as $(i,x)(i,y)=(i,xy)$ and $(i,x)^{-1}=(i,x^{-1})$ respectively, where $xy$ and $x^{-1}$ are the product  and inverse in the group $G_i$.
Then $\G$ becomes a groupoid with the unit space $\unsp=\{(i,e_{G_i}):\ i\in I\}$, where $e_{G_i}$ is the neutral element of $G_i$.

%\label{def:topgroupoiddef}
A \emph{topological groupoid} is a groupoid $\G$ endowed with a topology in which
the inverse map $\G\to \G$, $\gamma \mapsto \gamma^{-1}$ and the product map $\prodsp\to \G$, $(g,h)\to gh$ are continuous. Here, $\prodsp$ is equipped with the relative topology as a subset of $\G\times \G$. 
A \emph{locally compact} (respectively \emph{Hausdorff}) groupoid is a topological groupoid whose topology is locally compact (respectively Hausdorff).
Clearly, the range and source maps of a topological groupoid are continuous. 
If we take $\G$ to be Hausdorff, then $\unsp$ is also Hausdorff, and consequently, the fibers $\G^u$ and $\G_u$ are closed subsets of $\G$.

Locally compact groups are guaranteed Haar measures which allow a theory of integration to be developed. The analogous concept in the locally compact groupoid setting is called a {\it Haar system}. 
Roughly speaking, a Haar system is a family of positive, regular Borel measures defined on the fibers $\G^u$ satisfying a natural translation-invariance property. 
Unlike locally compact groups, however, locally compact groupoids do not necessarily admit Haar systems.
For a general discussion of Haar systems, see \cite[Definition 2.3]{DeG} or \cite[Definition 2.2.2]{PatersonBook}. 

Certain classes of groupoids, such as \'{e}tale groupoids, admit canonically defined Haar systems.
A locally compact, Hausdorff groupoid is called \emph{\'{e}tale} if the range map $r:\G\to \G$ (or equivalently, the source map) is a local homeomorphism.
All discrete groupoids are \'{e}tale, and if $\G$ is an \'{e}tale groupoid, each $\G^u$ is a discrete set in the relative topology.
All \'{e}tale groupoids have Haar systems consisting of counting measures, that is, each $\G^u$ is equipped with its counting measure. \'{E}tale groupoids have received a great deal of attention (see, for example, \cite{MR4533452,Brown-etal-2014,MR4592883,MR4510931,MR4645718,MR4608439,Matui-PAMS-2012,MR4483992,MR4563262,Pedro}). Unless otherwise specified, throughout this manuscript, we will assume our groupoids to be \'{e}tale.

\begin{rmk}
The terminology \emph{$r$-discrete} is not uniform in the literature and is sometimes used synonymously with \emph{\'{e}tale}. Throughout this paper, we use only the latter term; see \cite{ToolKit} for a discussion of the terminology.
\end{rmk}

%%%%%%%%%%%%%%%%%%%%%%%%%%%%%%%%%%%%%%%%%%%%%%%%%%%%%%%%%%%%%%%%%%%%%%%%%%%%%%%%%%%%%%%%%%%%%%%%
\subsection{Fourier and Fourier--Stieltjes algebras for groupoids}
\label{sec:groupoidalgebras}
To extend the definition of Fourier and Fourier--Stieltjes algebras associated with locally compact groups to the realm of groupoids, we start with defining what a representation of a locally compact groupoid is. 
Roughly speaking, while a group requires unitary representations over Hilbert spaces, groupoids will require $\G$-Hilbert bundles.

\begin{defn}[Continuous Banach field] \cite[Definition 10.1.2]{Dixmier}
\label{def:BanachField}
    A \emph{continuous field of Banach spaces} $\mathscr{E}$ over a topological space $T$ is a family $\set{E_t}_{t\in T}$ of Banach spaces, with a set $\Gamma \subseteq \Pi_{t\in T} E_t$ such that:
    \begin{enumerate}
        \item \label{linearspace} $\Gamma$ is a complex linear subspace of $\Pi_{t\in T} E_t$;
        \item \label{density} For every $t\in T$, the set $\set{x(t):x\in \Gamma}$ is dense in $E_t$;
        \item \label{normcontinuity} For every $x\in \Gamma$, the function $t\to \nnorm{x(t)}_{E_t}$ is continuous;
        \item \label{continuoussection} Let $x\in \Pi_{t\in T} E_t$.  If for every $t\in T$ and every $\epsilon>0$, there exists $x'\in \Gamma$ such that $\nnorm{x(\tau)-x'(\tau)}_{E_\tau}<\epsilon$ for all $\tau$ in some neighborhood of $t$, then $x\in \Gamma$.
    \end{enumerate}
    We refer to elements of $\Pi_{t\in T} E_t$ as \emph{sections}, and elements of $\Gamma$ as \emph{continuous sections}.
    In the case where each Banach space is also a Hilbert space, the definition above is referred to as a \emph{continuous Hilbert field}.
\end{defn}

In practice, continuous fields of Banach spaces are built as a `closure' of a collection of sections satisfying \eqref{linearspace}, \eqref{density} and \eqref{normcontinuity} of \Cref{def:BanachField}. This is done using the following result.

\begin{prop}\cite[Proposition 10.2.3]{Dixmier}\label{prop:fundamentalsections}
    Let $\set{E_t}_{t\in T}$ be a family of Banach spaces indexed over a topological space $T$. Suppose that there exists a set $\Xi \subseteq \prod_{t\in T}E_t$ of sections satisfying the first three conditions of \Cref{def:BanachField}. Then there is a unique subset $\Gamma\subseteq \prod_{t\in T}E_t$ that satisfies all four conditions and contains $\Xi$.

    The collection $\Xi$ is referred to as the \emph{fundamental family of sections}, and elements in $\Xi$ are called the \emph{fundamental sections}. 
\end{prop}

Moving forward, let $\fH:=(\set{\fH_t}_{t\in T},\Gamma)$ be a continuous Hilbert field over the space $T$. Denote the norm of each Hilbert space $\mathfrak{H}_t$ as $\|\cdot\|_t$ or $\|\cdot\|_{\mathfrak{H}_t}$ for $t\in T$. 
%The indices will be omitted when the norm is clear from the context. 
We define the \emph{section norm}
$$\nnorm{\xi}_\Delta=\sup_{t\in T}\nnorm{\xi(t)}_t.$$
Let $\Delta_b(\fH)$ denote the subspace of $\Gamma$ consisting of sections whose $\|\cdot\|_\Delta$-norm is bounded.  Equipped with pointwise operations and $\|\cdot\|_\Delta$-norm, the space $\Delta_b(\fH)$ is a Banach space (\cite[pg 6]{PatersonPaper}).
We further define  $\Delta_c(\fH)$ to be the subspaces of $\Gamma$ of continuous sections with compact support, and $\Delta_0(\fH)$ to be those sections in $\Gamma$ that vanish at infinity under the section norm. All three of $\Delta_b,\Delta_c$, and $\Delta_0$ are $C(T)$-modules under pointwise multiplication as the module action (\cite[Proposition 10.1.9]{Dixmier}).

We now give the definition of a groupoid representation as presented in \cite[Section 3]{PatersonPaper}.

\begin{defn}[Groupoid representation]\label{def:ghilbdle}
    Let $\G$ be a topological groupoid, and $\fH:=(\set{\fH_u}_{u\in \unsp},\Gamma)$ be a continuous Hilbert field. A \emph{(continuous unitary) representation of $\G$ on the continuous Hilbert field $\fH$} is a collection of unitary maps $\calpi(x):\fH_{s(x)}\to \fH_{r(x)}$, for each $x\in \G$,  that satisfy the following:
    \begin{enumerate}
        \item \label{coefficientcontinuity} For each $\bxi, \boldeta \in \Delta_b(\mathfrak{H})$ the map $\G\to \C$ defined as $x\mapsto \innprod{\calpi(x)\bxi(s(x)),\boldeta(r(x))}_{\fH_{r(x)}}$ is continuous;
        \item  For every $(x,y)\in \prodsp$, we have $\calpi(xy)=\calpi(x)\calpi(y)$. 
    \end{enumerate}
    A representation as above is also called a \emph{$\G$-Hilbert bundle over $\unsp$}. We use these two terminologies interchangeably. 
\end{defn}

The map $x\mapsto \innprod{\calpi(x)\bxi(s(x)),\boldeta(r(x))}_{\fH_{r(x)}}$ given in \Cref{def:ghilbdle} is referred to as the \emph{coefficient function of $\fH$ associated with $\bxi$ and $\boldeta$}, and denoted by $\fH_{\bxi,\boldeta}$.  To differentiate between coefficient functions of group representations, and coefficients of $\G$-Hilbert bundles, we will reserve lower case Greek letters and unbolded symbols for group representations, and accented text with bold symbols to refer to $\G$-Hilbert bundles and their sections. 
We note that it suffices to check condition \eqref{coefficientcontinuity} of \Cref{def:ghilbdle} for $\bxi,\boldeta\in \Xi$, the fundamental sections, that are $\|\cdot\|_\Delta$-bounded.

\begin{example}\label{exa:L2(G)}
       Let $\G$ be a locally compact \'{e}tale groupoid. Each $\G^u$ is a discrete space and is therefore naturally equipped with the counting measure.
       We define $L^2(\G)$ to be the continuous Hilbert field associated to the family of Hilbert spaces 
       $\set{\ell^2(\G^u)}_{u\in \unsp}$, with $C_c(\G)$ serving as a fundamental family of sections via the identification
\[
C_c(\G)\ni f \longmapsto {\{f^u\}}_{u\in \unsp},
\]
where $f^u=f|_{\G^u}$. 
It is straightforward to verify that $C_c(\G)$ indeed forms a fundamental family of sections in the sense of \Cref{prop:fundamentalsections}, and hence determines the continuous structure uniquely.
An equivalent description of $L^2(\G)$ is given in \cite[p.~10]{PatersonPaper}: a section $\{F^u\}_{u\in \unsp}$ is continuous if and only if the functions 
$u\mapsto \|{F^u}\|_{\ell^2(\G^u)}$ is continuous and $u\mapsto\innprod{F^u,g^u}$ is continuous for all $g\in C_c(\G)$. 
(See also \cite[Example 2.7]{DeG} for the definition of $L^2(\G)$ for general groupoids equipped with Haar systems.)

The \emph{left regular representation} of $\G$ is defined on the continuous Hilbert field $L^2(\G)$ as follows:
    For $x\in \G$, define 
    $$L_x: \ell^2(\G^{s(x)})\to \ell^2(\G^{r(x)}), \quad L_x\delta_z=\delta_{xz}.$$ 
    This will turn $L^2(\G)$ into a $\G$-Hilbert bundle, denoted $\Lambda$. 
\end{example}

\begin{defn}(\cite{PatersonPaper})\label{def:BG} %[pg. 7 lines 24-25, 30-32]
    The \emph{Fourier--Stieltjes algebra} of a locally compact groupoid $\G$ is defined as
    $$B(\G):=\set{\fH_{\bxi,\boldeta}:\fH \text{ is a $\G$-Hilbert Bundle, }\bxi,\boldeta\in \Delta_b(\fH)},$$
   equipped with pointwise operations and the norm defined as
   $$\nnorm{\varphi}_{B(\G)}:=\inf_{\varphi=\fH_{\bxi,\boldeta}}\nnorm{\bxi}_\Delta\cdot{\nnorm{\boldeta}_\Delta}$$
where the infimum is taken over all coefficient function representations of $\varphi$.
\end{defn}
It was shown in \cite{PatersonPaper} that $B(\G)$ is a unital commutative Banach algebra.
When $\G$ admits a Haar system, the Fourier algebra $A(\G)$ was also defined in \cite[page 15]{PatersonPaper}. 
Here, we recall this definition in the simpler setting of \'{e}tale groupoids.

\begin{defn}\label{def:FourierAlgebra}
Let $\G$ be a locally compact \'{e}tale groupoid, and $\Lambda$ denote the left regular representation of $\G$ on $L^2(\G)$.
Recall that $C_c(\G)$ is the set of fundamental sections of $L^2(\G)$.
Define
$$A_{cf}(\G)=\left\{\Lambda_{f,g}:\ f,g\in C_c(\G)\right\}.$$
Furthermore, let $A_c(\G)$ be the algebra generated by $A_{cf}(\G)$ under pointwise operations. The \emph{Fourier Algebra} of $\G$ is defined as $$A(\G):=\overline{A_c(\G)}^{\nnorm{\cdot}_{B(\G)}}.$$
\end{defn}
It is easy to see that both $A_{cf}(\G)$ and $A_c(\G)$ are contained in $C_c(\G)$. 
By definition, $A(\G)$ is a closed subalgebra of $B(\G)$ and is therefore itself a Banach algebra. 
As in the group case, $A(\G)\subseteq C_0(\G)$ by \cite[Proposition 14]{PatersonPaper}. Moreover, $A(\G)$ is an ideal of $B(\G)$ by \cite[Corollary 3]{PatersonPaper}, and \cite[Proposition 16(ii)]{PatersonPaper} shows that $A(\G)$ strongly separates points.
\begin{rmk}
There are at least two alternative proposed definitions of $A(\G)$. 
The first defines $A(\G)$ as the closure of the coefficient functions $\Lambda_{f,g}$ arising from sections $f,g$ that vanish at infinity. 
The second, due to \cite[page 186]{Oty}, defines $A(\G)$ as the closure of $C_c(\G)\cap B(\G)$ in the $B(\G)$-norm. It should be noted, however, that \cite{Oty} uses a different norm on the Fourier--Stieltjes algebra from the one used here.  
\end{rmk}

Next, we show that the two proposed definitions of $A(\G)$ inspired by the above remark coincide for \'{e}tale groupoids. 
\footnote{This answers a question posed by Brian Forrest during the Ph.D. defense of the first-named author.}
\begin{lem}\label{lem-etal-Cc-Acf} 
    Let $\G$ be a (locally compact Hausdorff) \'{e}tale groupoid. Then, $A_{cf}(\G)=A_c(\G)=C_c(\G)$ and 
    $A(\G)=\overline{B(\G)\cap C_c(\G)}^{\nnorm{\cdot}_{B(\G)}}$.
\end{lem}
\begin{proof}
By definition, $A_{cf}(\G)\subseteq A_c(\G)\subseteq C_c(\G)$. 
To prove the reverse inclusions, let $f\in C_c(\G)$ be given, and note that $s(\operatorname{supp}f)$ is a compact subset of $\G^{(0)}$, as the source map $s$ is continuous.
Choose a real-valued $h\in C_c(\G^{(0)})$ equal to $1$ on $s(\operatorname{supp}f)$, and extend $h$ to $\widetilde{h}:\G\to {\mathbb C}$ as follows:
$$\widetilde{h}(x)=\left\{\begin{array}{ll}
    h(x) &  x\in \G^{(0)}\\
    0 & \mbox{otherwise} 
\end{array}\right.$$
Note that $\G$ is Hausdorff and \'{e}tale, so the unit space $\G^{(0)}$ is an open and closed subset of $\G$. Therefore, the function $\widetilde{h}$ belongs to $C_c(\G)$.
Moreover, for every $\gamma\in \G$, we have
$$
\Lambda_{\widetilde{h},\overline f}(\gamma)
=
\sum_{t\in\G^{r(\gamma)}}
\widetilde{h}(\gamma^{-1}t)f(t)
=h(\gamma^{-1}\gamma)f(\gamma)
=h(s(\gamma))f(\gamma)
=f(\gamma),
$$
where the second equality follows from the fact that $t=\gamma$ is the only possible element of $\G$ that satisfies
$\gamma^{-1}t\in \supp(h)\subseteq \G^{(0)}$.
Since both $\widetilde{h}$ and $\overline{f}$ are elements of $C_c(\G)$, we get
$\Lambda_{\widetilde{h},\overline f}\in A_{cf}(\G)$. This establishes $A_{cf}(\G)=A_c(\G)=C_c(\G)$.
Finally, we conclude 
    $$\overline{B(\G)\cap C_c(\G)}^{\nnorm{\cdot}_{B(\G)}}=\overline{B(\G)\cap A_c(\G)}^{\nnorm{\cdot}_{B(\G)}}=\overline{A_c(\G)}^{\nnorm{\cdot}_{B(\G)}}=A(\G).$$
\end{proof}

%%%%%%

\begin{defn}[Restriction maps]\label{def:restriction-map}
    For a locally compact groupoid $\G$ and an element $u\in \unsp$, we define the map $R^u_{B(\G)}: B(\G)\to B(\G_u^u)$ by $\varphi \mapsto \varphi\mid_{\G_u^u}$. An analogous restriction function is defined on $A(\G)$ via $R^u_{A(\G)}: A(\G)\to A(\G_u^u)$, $\varphi \mapsto \varphi\mid_{\G_u^u}$. 
\end{defn}

The above restriction maps were introduced in \cite{DeG}. 
Clearly, $R^u_{A(\G)}=R^u_{B(\G)}\mid_{A(\G)}$. 
However, the two maps exhibit different fundamental properties. 
For example, when $\G$ is \'{e}tale, $R^u_{A(\G)}$ is surjective, while $R^u_{B(\G)}$ need not be. 
We therefore used the subscripts to distinguish their domains. In this paper, we consider only restrictions of elements of $A(\G)$ and accordingly simplify our notation to $R_u:A(\G)\to A(\G_u^u)$.

%%%%%%%%%%%%%%%%%%%%%%%%%%%%%%%%%%%%%%%%%%%%%%%%%%%%%%%%%%%%%%%%%%%%%%%%%%%%%%%%%%%%%%%%%%%%%%%%%%%%%%%%%%%%%%%%%%%%%%%%%%%%%%%%%%%%%%%%%%%%%%%%%%%%%%%%%%
\section{HLS groupoids and their Fourier algebras}\label{sec:HLSAlgebras}

In this section, we study a particular class of group bundles called the HLS\footnote{HLS stands for Higson, Lafforgue, and Skandalis, the authors who originated the construction.} groupoids, originally appearing in \cite{Skandalis}. 
The definition of HLS groupoids is based on the notion of an approximating sequence of normal subgroups, which we now recall.
An \emph{approximating sequence} for a discrete group $G$ is a sequence $\set{K_n}_{n\in \N}$ of normal subgroups of $G$ with finite index satisfying 
$$\bigcap_{n\in \N} K_n=\set{e_G}\ \mbox{ and }\ K_n \supseteq K_{n+1}\ \mbox{ for } n\in \N.$$
The following lemma guarantees the existence of approximating sequence of normal subgroups for a large class of groups. This lemma is a generalization of \cite[Lemma 2.8]{Willett}, where a similar result for the group $\F_2$ is proved. Recall that a \emph{residually finite} group is one in which if $x\neq e_G$ then there exists a homomorphism $\varphi$ from $G$ to a finite group such that $\varphi(x)\neq e$.
\begin{lem}\label{lem:apxsequence}
    Let $G$ be a residually finite and finitely generated group. Then $G$ has an approximating sequence.
\end{lem}
\begin{proof}
    Since $G$ is finitely generated, there are only finitely many homomorphisms from $G$ into groups of order at most $n$. Define $K_n$ to be the intersection of the kernels of all such homomorphisms. Each $K_n$ is then a finite-index normal subgroup of $G$, because finite intersections of finite-index subgroups again have finite index. The sequence $\{K_n\}$ is nested by construction, and residual finiteness of $G$ implies that the intersection of all the sets $K_n$, for $n\in {\mathbb N}$ is trivial.
Hence $\{K_n\}$ is an approximating sequence for $G$.
\end{proof}

\begin{defn}
    Given a discrete group $G$ and an approximating sequence $\{K_n\}_{n\in \N}$, let $G_n$ denote the finite quotient group $G/K_n$, and  $q_n: G\to G_n$ be the canonical quotient map. We set $G_\infty:=G$ and $q_\infty$ the identity map. 
The \emph{HLS groupoid} associated to $(G,\set{K_n}_{n\in \N})$ is the group bundle
    $$\G:=\sqcup_{n\in \N\cup\{\infty\}}\set{n}\times G_n,$$
where $\G$ is equipped with the topology generated by the following open sets:
\begin{enumerate}
    \item[(i)] The singleton $\set{(n,g)}$ for each $n\in \N$ and $g\in G_n$;
    \item[(ii)] The \emph{ray} at $(\infty,g)$ defined as  $\ray_n(g):=\set{(m,q_m(g)): m\in \N\cup\{\infty\}, m\geq n}$ for fixed $g\in G$ and $n\in \N$.
\end{enumerate}
\end{defn}

As any HLS groupoid $G$ is a group bundle, we will often think of each individual group $G_n$ as a subset of $\G$ via the identification $q_n(x)\sim (n,q_n(x))$. 
HLS groupoids are locally compact, Hausdorff and \'{e}tale. 
If the discrete group $G$ is countable, then the associated HLS groupoid $\G$ is second countable. 
The topology of an HLS groupoid gives a particularly simple description of continuous functions. Indeed, a function $f:\G\to\C$ is continuous if and only if it is continuous along each ray, i.e., for every $g\in G_\infty$, one has
\begin{equation}\label{eq:cts-function}
    \lim_{n\to\infty} f(n,q_n(g)) = f(\infty,g).
\end{equation}

The following result collects useful information about the topology of an HLS groupoid and the description of its compact subsets.
\begin{prop}\label{lem:rays}
 Let $\G$ be the HLS groupoid associated to the pair $(G,\set{K_n}_{n\in \N})$. Then we have the following properties:
 \begin{enumerate}
     \item\label{lem-item:ray1} For every $n\in \N$ and $g\in G$, the ray $\ray_{n}(g)$ is an open and closed compact set.
     \item\label{lem-item:ray3} A set $C\subseteq \G$ is compact if and only if there exists a finite subset $\{g_1,\ldots,g_n\}$ of $G$ and a finite subset $F$ of $C$ such that %\label{prop:HLSTopology}
    \begin{equation}\label{eq:cpt-HLS}
        \{(\infty,g_1),\ldots,(\infty,g_n)\}\subseteq C\setminus F\subseteq \bigcup_{1\leq i\leq n}\ray_1(g_i).
    \end{equation}
     \item\label{lem-item:ray2} Let $E$ be a finite set in $G$. There exists $N\in \N$ such that $\ray_N(x)\cap \ray_N(y)=\varnothing$ for every $x\neq y\in E$.
 \end{enumerate}
\end{prop}
\begin{proof}
    Let $\set{U_\alpha}_{\alpha\in I}$ be an open cover of $\ray_n(x)$. At least one set in the open cover, say $U_{\alpha_x}$, must contain $(\infty,x)$. By the definition of the topology of $\G$, there exists $n'\in \N$ such that $\ray_{n'}(x)\subseteq U_{\alpha_x}$. The set $\ray_{n}(x)\setminus \ray_{n'}(x)$ is finite (or possibly even empty), and can be covered by a finite number of sets in the open cover $\set{U_\alpha}_{\alpha\in I}$. This proves the compactness of the ray.
    As $\G$ is Hausdorff and each $\ray_n(x)$ is compact, they are also closed. 
    This finishes the proof of Item \Cref{lem-item:ray1}.

    To prove Item \ref{lem-item:ray3}, let $C$ be a subset of $\G$ and consider the following open cover of $C$:
$$\Big\{V_g:=\ray_1(g):\ (\infty,g)\in C\Big\}\bigcup\Big\{W_g:=\{(n,g)\}:\ n\in\N,\ g\in G_n,\ (n,g)\in C \Big\}.$$
Any subcover of the above open cover should contain all sets of the form $V_g$. So, if $C$ is compact, then $C\cap G_\infty$ must be finite. In addition, to find a finite subcover of the above cover, $C\setminus (\cup_{(\infty,g)\in C\cap \G_\infty}V_g)$ must be finite as well. This proves the `only if' direction. 

To prove the `if' direction of the statement, fix an element $(\infty, g)$ in $\G$ with $g\in G$. Note that every open set in $\G$ containing $(\infty, g)$ should include some ray at $(\infty, g)$ as a subset. 
Suppose $C$ is a subset of $\G$ satisfying \eqref{eq:cpt-HLS}, and let $\{O_\alpha:\ \alpha\in I\}$ be an open cover for $C$.
Let $1\leq i\leq n$.
From the previous discussion, since $(\infty,g_i)\in C$, there exists $k_i\in \N$ and $\alpha_i\in I$, such that $\ray_{k_i}(g_i)\subseteq O_{\alpha_i}$. Thus, $\{O_\alpha:\ \alpha\in I\}$ has a finite subcover, and $C$ is compact.

To prove Item \ref{lem-item:ray2}, note that the set $\set{x^{-1}y:x,y\in E, x\neq y}$ is a finite set not including $e_{G_\infty}$. Since $\cap_{n\in \N} K_n=\set{e_{G_{\infty}}}$,  for every $z\neq e_{G_\infty}$ there exists an index $N_z$ such that $z\notin K_{N_z}$. Let 
    $$N=\max\set{N_{x^{-1}y}: x,y\in E, \ x\neq y}.$$ 
    Since $\set{K_n}_{n\in \N}$ is a nested sequence of normal subgroups, for each $n\geq N$, we have $x^{-1}y\notin K_n$, and consequently, $q_n(x)\neq q_n(y)$. Therefore, $\ray_{N}(x)\cap \ray_N(y)=\emptyset$.
\end{proof}

\subsection{Fourier algebra of an HLS groupoid}\label{subsec:fourier-hls} 
Here, we present the simplified definitions related to the representations and Fourier algebra of HLS groupoids.
To begin, we define $\infN:=\N\cup\set{\infty}$ and equip it with the topology of one point compactification of $\N$.

\begin{defn}[\Cref{def:BanachField} simplified]
\label{def:HLSHilbertField}
    A {\it continuous Hilbert field} indexed over $\infN$ is a family of Hilbert spaces $\{\Hil_n\}_{n\in \infN}$ together with a set $\Gamma\subseteq \Pi_{n\in \infN}\Hil_n$ satisfying the following properties: 
    \begin{enumerate}
        \item[$1'$.] $\Gamma$ is a complex linear subspace of $\Pi_{n\in \infN}\Hil_n$;
        \item[$2'$.] For every $n\in \infN$, the set $\set{x(n):x\in \Gamma}$ is dense in $\Hil_n$;
        \item[$3'$.] For every $x\in \Gamma$, $\lim_{n\to \infty}\nnorm{x(n)}_{\Hil_n}=\nnorm{x(\infty)}_{\Hil_\infty}$;
        \item[$4'$.] Let $x\in \Pi_{n\in \infN}\Hil_n$. If, for every $\epsilon>0$, there exists an $x'\in \Gamma$ and an $N\in \N$ such that $\nnorm{x'(n)-x(n)}_{\Hil_n}<\epsilon$ for every $n\in [N,\infty]$ then $x\in \Gamma$.
        \item[$5'$.] For every $n\in \N$, we have $\set{x(n):x\in \Gamma}=\Hil_n$. In particular, $\Hil_n$ embeds naturally in $\Gamma$ in the following sense: 
        $$\widetilde{\Hil_n}:=\left\{\{x(k)\}_{k\in \N^*}: x(n)\in \Hil_n \mbox{ and } x(k)=0 \ \forall k\neq n\right\}\subseteq \Gamma.$$
    \end{enumerate}
\end{defn}
Note that the above definition is equivalent to \Cref{def:BanachField} for an HLS groupoid. Indeed, since $\N$ is a discrete subset of $\infN$, the density condition given in Item $2'$ implies that Condition 4 of \Cref{def:BanachField} is equivalent to Items $4'$ and $5'$.
Moreover, for every $\bxi\in \Gamma$, the supremum
$\nnorm{\bxi}_\Delta=\sup_{n\in \infN}\nnorm{\bxi(n)}_{\Hil_n}$
is bounded due to Item $3'$. This demonstrates that the set of all bounded sections matches the set of continuous sections, i.e., $\Delta_b=\Gamma$.

\begin{defn}[\Cref{def:ghilbdle} simplified]\label{def:HLSRepresentation}
    Let $\G$ be an HLS groupoid. 
    A {\it continuous unitary representation} $\varpi$ of $\G$ consists of a continuous Hilbert field $(\set{\Hil_n}_{n\in \infN},\Gamma)$ together with, for each $n\in \infN$, a unitary representation $\pi_n:\G\to {\cal U}(\Hil_n)$ 
    such that the associated matrix coefficient functions vary continuously across the field. 
    More precisely, for every $\bxi,\boldeta\in \Gamma$, the map
    $$\varpi_{\bxi,\boldeta}:\G\to \C, \;\;\; \varpi_{\bxi,\boldeta}(n,q_n(x))=\innprod{\pi_n(q_n(x))\bxi(n),\boldeta(n)}_{\Hil_n}$$
    is continuous.
\end{defn}
Using \Cref{eq:cts-function}, the above coefficient function is continuous precisely when for each $x\in G_\infty$, we have
$$\lim_{n\to \infty}\innprod{\pi_n(q_n(x))\bxi(n),\boldeta(n)}_{\Hil_n}=\innprod{\pi_\infty(x)\bxi(\infty),\boldeta(\infty)}_{\Hil_\infty}.$$

Let $\Lambda$ be the left regular representation of the HLS groupoid $\G$ with the collection of continuous sections $\Gamma\subseteq \prod_{n\in\infN}\ell^2(G_n)$.
Recall that $\Gamma$ is formed via the unique closure of $C_c(\G)$ (see \cite[Proposition 10.2.3]{Dixmier} and \Cref{exa:L2(G)}).
Every continuous section $\{f_n\}_{n\in \infN}$ can be naturally identified with a function $F:\G\to \C$ as follows:
\begin{equation}\label{eq:section-to-function}
    F(n,q_n(g))=f_n(q_n(g)), \mbox{ for }\ g\in G_\infty.
\end{equation}

We now provide a concrete description of the continuous sections of the left regular representation for the case of HLS groupoids, and show that the coefficient functions associated with any pair of continuous sections in $\Gamma$ belong to the groupoid Fourier algebra. 
\begin{prop}\label{cor-cts-sec}
    Let $\G$ be an HLS groupoid, and $\Gamma$ the set of continuous sections of the left regular representation  of $\G$. 
    Consider the natural identification of sections in $\Gamma$ with functions on $\G$ given in \eqref{eq:section-to-function}.
    Then, we have the following.
    \begin{enumerate}
        \item\label{item:section1} For every $\{f_n\}_{n\in \infN}\in \Gamma$, the associated function $F$ is continuous on $\G$.
        \item\label{item:section2} $\Gamma=\left\{F\in C(\G):\ F|_{G_\infty}\in \ell^2(G_\infty), \ \nnorm{F|_{G_n}}_2\to \nnorm{F|_{G
        _\infty}}_2\right\}.$
        Here we have identified $G_n$, for $n\in \N^*$, with its image $\{n\}\times G_n$ in $\G$.
        \item\label{item:sections-coefficients} 
For every $F,H\in\Gamma$, we have $\Lambda_{F,H}\in A(\G)$.
    \end{enumerate}
\end{prop}

\begin{proof}
    Let $\{f_n\}_{n\in \infN}\in \Gamma$ be given, and fix $\epsilon>0$ and $g\in G_\infty$.
    As described in \Cref{exa:L2(G)}, $\Gamma$ is the closure of the space of the fundamental sections $C_c(\G)$. 
    Hence, there exists $h\in C_c(\G)$ and $N_1\in \N$ such that
    $$\nnorm{f_m-h_m}_2<\frac{\epsilon}{3}\ \mbox{ for }\ m\in [N_1,\infty],$$
    where $h_m=h|_{G_m}$.
    Consequently, for every $m\in [N_1,\infty]$ and every $g\in G_\infty$, we have 
    $$\abs{f(m,q_m(g))-h(m,q_m(g))}<\frac{\epsilon}{3}.$$ 
    Since $h\in C_c(\G)$, there exists $N_2\in \N$ such that $\abs{h(n,q_n(g))-h(\infty,g)}<\frac{\epsilon}{3}$ for each $n\geq N_2$. Applying an $\epsilon/3$-argument, for $N\geq \max\set{N_1,N_2}$, 
    we have $\abs{f(n,q_n(g))-f(\infty,g)}\leq \epsilon$. This proves \Cref{item:section1}.

    To prove \Cref{item:section2}, note that the inclusion `$\subseteq$' follows directly from part \Cref{item:section1} and the definition of a continuous Hilbert field.
    Conversely, suppose that $F:\G\to \C$ satisfies $F\in C(\G)$,  $F|_{G_\infty}\in \ell^2(G_\infty)$, and 
    $\nnorm{F|_{G_n}}_2\to \nnorm{F_{G_\infty}}_2$. 
    By the criterion recalled in \Cref{exa:L2(G)}, it suffices to show that, for every fundamental section $H\in C_c(\mathcal G)$, we have
    $$ \langle F|_{G_n},H|_{G_n}\rangle_{\ell^2(G_n)} \longrightarrow \langle F|_{G_\infty},H|_{G_\infty}\rangle_{\ell^2(G_\infty)}. $$
    Since $H$ is compactly supported, there exists a finite set of elements of $G_\infty$, say $z_1,\ldots, z_m$, and a sufficiently large $N$, such that 
    \begin{equation}\label{eq:support-H}
        \supp(H\mid_{\G\setminus \sqcup_{i=1}^NG_i})\subseteq \sqcup_{i=1}^m \ray_N(z_i).
    \end{equation}
    So, by continuity of $F$ and $H$ we get:
    \begin{align*}
        \langle F|_{G_\infty},H|_{G_\infty}\rangle_{\ell^2(G_\infty)}
        &=\sum_{i=1}^m F(\infty,z_i)\overline{H(\infty, z_i)}\\
        &=\sum_{i=1}^m \big(\lim_{n\to \infty}F(n,q_n(z_i))\big)\big(\overline{\lim_{n\to \infty}H(n, q_n(z_i))\big)}\\
        &=\lim_{n\to \infty} \sum_{i=1}^m F(n,q_n(z_i))\overline{H(n, q_n(z_i))}\\
        &=\lim_{n\to \infty} \langle F|_{G_n},H|_{G_n}\rangle_{\ell^2(G_n)}, 
    \end{align*}
    where the last equality follows from \Cref{eq:support-H}.

    To prove \Cref{item:sections-coefficients}, we first show that $C_c(\G)$ is $\|\cdot\|_\Delta$-dense in $\Gamma$. Let $F\in\Gamma$ and $\epsilon>0$. Since $\Gamma$ is the (local) closure of the fundamental sections $C_c(\mathcal{G})$, there exist $K\in C_c(\G)$ and $N\in\N$ such that $\|F|_{G_n}-K|_{G_n}\|_2<\epsilon$ for all $n\in[N,\infty]$. Since $\sqcup_{n<N}G_n$ is a finite set, the function
    \[ K'(n,x):= 
    \begin{cases}
    F(n,x), & n<N,\\
    K(n,x), & N\leq n\le \infty,
    \end{cases}\] 
    belongs to $C_c(\G)$, and $\|F-K'\|_\Delta =\sup_{n\ge N}\|F|_{G_n}-K|_{G_n}\|_2\le\epsilon$. 
    This proves our claim. 
    
    Now let $\epsilon>0$ be arbitrary, and choose $F',H'\in C_c(\mathcal{G})$ with $\|F-F'\|_\Delta< \epsilon$ and $\|H-H'\|_\Delta<\epsilon$.  
    Clearly, $\Lambda_{F',H'}\in A_c(\G)$.
    Also, by the discussion after \Cref{def:HLSHilbertField}, we have  $\Gamma=\Delta_b$ and $\Lambda_{F,H}\in B(\G)$. The definition of the Fourier-Stieltjes norm gives that
    \[\|\Lambda_{F,H}-\Lambda_{F',H'}\|_{B(\G)} 
    \leq\|F-F'\|_\Delta\|H\|_\Delta+\|F'\|_\Delta\|H-H'\|_\Delta\leq \epsilon(\|H\|_\Delta+\|F\|_\Delta+\epsilon),\]
    which implies that $\Lambda_{F,H}\in \overline{A_c(\G)}^{\|\cdot\|_{B(\G)}}=A(\G)$.
\end{proof}

    Recall that for any \'{e}tale groupoid $\G$, we have $A_{c}(\G)=C_c(\G)$ and 
    $A(\G)=\overline{B(\G)\cap C_c(\G)}^{\nnorm{\cdot}_{B(\G)}}$ (\Cref{lem-etal-Cc-Acf}).
    When $\G$ is an HLS groupoid, the representation of $f\in C_c(\G)$ as a coefficient function in $A_c(\G)$ takes a simpler form. 
    For example, consider the case where $f\in C_c(\G)$ satisfies $\supp(f)\subseteq \mathds{1}_{\ray_n(x)}$ for some $x\in G_\infty$ and $n\in {\mathbb N}$.
    By \Cref{lem:rays} Item \ref{lem-item:ray1}, the function $g=\mathds{1}_{\ray_n(e_{G_\infty})}$ belongs to $C_c(\G)$.
    Then, for $m\in \infN$ and $z\in G_\infty$, we have
        $$\Lambda_{g,\overline{f}}(m,q_m(z))
        =\sum_{q_m(y)\in G_m} g(m,q_m(z^{-1}y))f(m,q_{m}(y))
        =\left\{\begin{array}{cc}
        f(m,q_{m}(z)) &  \mbox{ if } m\geq n\\
        0 & \mbox{ if } m<n    
        \end{array}\right.=f(m,q_m(z)).$$
Letting $f=\mathds{1}_{\ray_n(e_{G_\infty})}$ in the above discussion, we get the following easy corollary.
\begin{cor}\label{cor:characteristic-in-A}
   For every $x\in G_\infty$ and $n\in \N$, the characteristic function $\mathds{1}_{\ray_n(e_{G_\infty})}$ belongs to $A(\G)$. In fact, we have $\Lambda_{\mathds{1}_{\ray_n(e_{G_\infty})},\mathds{1}_{\ray_n(e_{G_\infty})}}=\mathds{1}_{\ray_n(e_{G_\infty})}.$
\end{cor}

\begin{rmk}\label{rmk:SubgroupoidRestriction}
Note that restriction-type results depend on the subobject being a group rather than a groupoid.
In particular, our results cannot be extended to restriction maps onto arbitrary subgroupoids of a given groupoid.
Indeed, let $\G$ be an HLS groupoid, and set
$\G_{\mathbb N}:=\G\setminus G_\infty.$
By \Cref{cor:characteristic-in-A}, for $x\in G$ and $n\in \N$, we have $\mathds{1}_{\ray_n(x)}\in A(\G).$
However,
\[
\mathds{1}_{\ray_n(x)}|_{\G_{\mathbb N}}\notin A(\G_{\mathbb N}).
\]
To see this, note that $\G_{\mathbb N}$ is an \'{e}tale groupoid, and hence
$A(\G_{\mathbb N})\subseteq C_0(\G_{\mathbb N}).$
Consequently, the restriction map $A(\G)\to  A(\G_{\mathbb N})$
is not well defined.
\end{rmk}

We end this section by establishing several results on the norms of the restricted functions.
\begin{prop}\label{prop:normconvergence1}
    Let $\G=\sqcup_{n\in \N^*} G_n$ be an HLS groupoid, and for each $n\in \infN$, let  $R_n:A(\G)\to A(G_n)$ denote the restriction map. 
    Let $u\in A(\G)$.
    \begin{enumerate}
        \item\label{item:lower-bnd-gen} $\sup\limits_{n\in \infN}\nnorm{R_n(u)}_{A(G_n)}\leq \nnorm{u}_{A(\G)}.$
        \item\label{item:upper-bnd-A-cf} $\|R_\infty(u)\|_{A(G_\infty)}\leq \liminf_{n}\|R_n(u)\|_{A(G_n)}.$
        \item\label{item:when-u-zero} Suppose $u\in A(\G)\cap C_c(\G)$. If $R_\infty(u)=0$  then 
        $\lim_{n\to \infty}\|R_n(u)\|_{A(G_n)}=0.$
    \end{enumerate}
\end{prop}
\begin{proof}
    Item~\ref{item:lower-bnd-gen} follows from \cite[Proposition 3.4(iii)]{DeG}, as each restriction map $R_n$ is contractive.
    To prove Item \ref{item:upper-bnd-A-cf}, we view each $R_n(u)$ as a function on $G_\infty$. 
    Namely, for $n\in \N$, define 
    $$\tilde{u}_n:G_\infty\to \C,\quad \tilde{u}_n(x)=u(n,q_n(x))=R_n(u) (q_n(x)).$$
    \begin{claim}\label{claim:fiber-in-infty}
        We have $\tilde{u}_n\in B(G_\infty)$ and $\|\tilde{u}_n\|_{B(G_\infty)}\leq \|R_n(u)\|_{A(G_n)}$.
        Consequently, $\{\tilde{u}_n\}_{n\geq N}$ is a bounded sequence in $B(G_\infty)$.
    \end{claim}
    \noindent\emph{Proof of claim.}
    Let $R_n(u)=\langle\lambda_{G_n}(\cdot)\xi_n,\eta_n\rangle$ be a representation of $R_n(u)$ as a coefficient function of the left regular representation of $G_n$ with $\xi_n,\eta_n\in \ell^2(G_n)$ and $\|R_n(u)\|_{A(G_n)}=\|\xi\|_2\|\eta_n\|_2$.
    It is easy to see that $\lambda_{G_n}\circ q_n$ is a unitary representation of $G_\infty$, 
    and for every $x\in G_\infty$,
    $\tilde{u}_n(x)=\langle(\lambda_{G_n}\circ q_n)(x)\xi_n,\eta_n\rangle$.
    So  $\tilde{u}_n\in B(G_\infty)$.
    Moreover,
    $$\|\tilde{u}_n\|_{B(G_\infty)}\leq \|\xi_n\|_2\|\eta_n\|_2=\|R_n(u)\|_{A(G_n)}.$$
    This finishes the proof of the claim, since we have already proved that the sequence $\{\|R_n(u)\|_{A(G_n)}\}_{n}$ is bounded above.
    
    As $u$ is a continuous function on $\G$, using \eqref{eq:cts-function}, we get $\tilde{u}_n\to R_\infty(u)$ pointwise. 
    By Eymard’s bounded convergence theorem for $B(G_\infty)$, the sequence $\{\tilde{u}_n\}_{n\in \N}$ converges to $R_\infty(u)$ in the weak* topology induced by $C^*(G_\infty)$. Since the norm is weak* lower semi-continuous, we have 
    \begin{equation}
        \|R_\infty(u)\|_{A(G_\infty)}\leq \liminf_{n}\|\tilde{u}_n\|_{B(G_\infty)} 
    \leq \liminf_{n}\|R_n(u)\|_{A(G_n)},
    \end{equation}
    where we used \Cref{claim:fiber-in-infty} in the last inequality. This proves Item \Cref{item:upper-bnd-A-cf}.

    To prove Item \ref{item:when-u-zero}, note that since $u\in C_c(\G)$, there are finitely many points of the form $(\infty, x_1),\dots, (\infty, x_m)$ contained in $\supp(u)$.
    Moreover, applying \Cref{lem:rays}, Items \Cref{lem-item:ray3} and \Cref{lem-item:ray2}, for sufficiently large $N$, we have 
    \begin{equation}\label{eq:finite-rays}
        \supp(u\mid_{\G\setminus \sqcup_{i=1}^NG_i})\subseteq \sqcup_{i=1}^m \ray_N(x_i),
    \end{equation}
    where $\sqcup$ denotes the union of disjoint sets.
    So for $n\in [N,\infty]$, each function $R_n(u)$ is supported in $\{(n,q_n(x_1)),\ldots,(n,q_n(x_m))\}$.
    On the other hand, since $u$ is continuous on $\G$, for every $x\in G_\infty$ we have
    $$u(n,q_n(x))\to u(\infty,x) \mbox{ as } n\to \infty.$$
    The above convergence holds simultaneously for $x_1,\ldots,x_m$. So, $\|R_n(u)\|_\infty\to 0$ as $n\to \infty.$
    Thus, for $n\geq N$ we have
    $$\|R_n(u)\|_{A(G_n)}=\|\sum_{i=1}^mu(n,q_n(x_i))\delta_{q_n(x_i)}\|_{A(G_n)}\leq \sum_{i=1}^m 
    |u(n,q_n(x_i))|\leq m\|R_n(u)\|_{\sup},$$
    where we used the fact that the Fourier norm of a point-mass function  is bounded above by 1.
    So, $\|R_n(u)\|_{A(G_n)}\to 0$ as desired. 
\end{proof}

%%%%%%%%%%%%%%%%%%%%%%%%%%%%%%%%%%%%%%%%%%%%%%%%%%%%%%%%%%%%%%%%%%%%%%%%%%%%%%%%
\section{Functorial properties of the Fourier algebras of HLS groupoids}\label{sec:FunctorialProperties}
In this section, we exploit the particular structure of HLS groupoids to strengthen general functorial results for the Fourier algebra of \'{e}tale groupoids in the setting of HLS groupoids. Namely, for each $n\in \N^*$, we construct explicit norm-preserving lifts for the restriction map $R_n:A(\G)\to A(G_n)$ (\Cref{thm:HLSSurjective}). 
%Second, using these lifts, we identify of the kernel of the restriction map $R_\infty$ with the $c_0$-direct sum of group Fourier algebras $A(G_n)$ (\Cref{cor:norm-convergence-fibers-kernel}).
%Finally, we extend the convergence result of \Cref{prop:normconvergence1} to the general setting (\Cref{cor:norm-A(G)-limit-sections}). 

When $n\in \N$, there is a natural embedding from each $A(G_n)$ into $A(\G)$ which also maintains the norm.
\begin{lem}\label{lem:spectrumembedding}
Let $n\in \N$.
    Define $\iota_n:A(G_n)\to A(\G)$ by 
    $$(\iota_nu)(m,q_m(x)):=\begin{cases}
        u(q_n(x)) & m=n\\
        0 & m\neq n
    \end{cases}.$$
    Then, $\iota_n$ is an algebra homomorphism, which is also an isometry.
    Moreover, $\iota_n(A(G_n))$ is an ideal in $A(\G)$.
\end{lem}
\begin{proof}
    Clearly $\iota_n$ is both linear and multiplicative. If $u\in A(G_n)$, then we may choose $f_1,f_2\in \ell^2(G_n)$ such that $\nnorm{u}_{A(G_n)}=\nnorm{f_1}_2\nnorm{f_2}_2$. Moreover, by defining $F_i:\G\to \C$ via $F_i(n,q_n(x))=f_i(q_n(x))$ and $F_i(m,q_m(x))=0$ if $m\neq n$, we see that $F_i\in C_c(\G)$ and thus defines continuous and bounded sections of $\Lambda$. The norm $\nnorm{F_i}_\Delta$ is precisely $\nnorm{f_i}_2$, and so we have
    $$\nnorm{\iota_n(u)}_{A(\G)}=\nnorm{\Lambda_{F_1,F_2}}_{A(\G)}\leq \nnorm{F_1}_\Delta \nnorm{F_2}_\Delta=\nnorm{f_1}_2\nnorm{f_2}_2=\nnorm{u}_{A(G_n)}.$$
    By \cite[Proposition 3.4]{DeG}, the restriction map $R_n$ is contractive, and so 
    $$\nnorm{u}_{A(G_n)}=\nnorm{R_n(\iota_n(u))}_{A(G_n)}\leq \nnorm{\iota_n(u)}_{A(\G)}$$
    demonstrating $\iota_n$ is an isometry. 

    Let $u\in A(G_n)$ and $v\in A(\G)$ be arbitrary. We want to show that $\iota_n(u)v\in \iota_n(A(G_n))$.
    For every $x\in G_\infty$, we have
    $$(\iota_n(u)v)(m,q_m(x))=\begin{cases}
        0 & m\neq n\\
        u(q_n(x))v(n,q_n(x)) & m=n
    \end{cases},$$
    which gives $\iota_n(u)v=\iota_n(uR_n(v))$, which finishes the proof. 
\end{proof}
In view of the above lemma, we will regard each $\iota_n$ as an embedding of $A(G_n)$ into $A(\G)$.
\begin{thm}\label{thm:HLSSurjective}
    Let $\G$ be the HLS groupoid associated to the pair $(G,\set{K_n}_{n\in \N})$. Then for each $n\in \N\cup\{\infty\}$ and $u\in A(G_n)$, there exists $v\in A(\mathcal G)$ such that $R_n(v)=u$ and $\|v\|_{A(\mathcal G)}=\|u\|_{A(G_n)}$. In particular, the map $R_n:A(\G)\to A(G_n)$ is a quotient map.
\end{thm}

\begin{proof}
    Recall that we set $G_\infty:=G$, and denote by $q_n:G \to G/K_n$ the canonical quotient map. 
    If $n<\infty$, the restriction map $R_n$ is a quotient map, since any $f\in \ell^2(G_n)$ embeds isometrically as an element of $C_c(\G)\subseteq \Gamma$ as described in the proof of \Cref{lem:spectrumembedding}. So, we will only focus on the map $R_\infty$ and show that every $L_{f,g}\in A(G_\infty)$ with $f,g\in \ell^2(G_\infty)$ has a pre-image under $R_\infty$ whose norm is equal to $\|L_{f,g}\|_{A(G_\infty)}$.

    To begin, fix $f\in \ell^2(G_\infty)$, and assume that $f\geq 0$ pointwise. For $n\in \N$ define the following functions:
    \begin{align}\label{eq:HLSSection1}
        &f_n :G/K_n\to \C\\   
        &f_n(q_n(x)) :=\nnorm{f|_{xK_n}}_2, \quad x\in G \nonumber 
    \end{align}
    Next define
    \begin{align}\label{eq:HLSSection2}
        &F:\G \to \C \\
        &F(n,q_n(x)):= \begin{cases}
           f_n(q_n(x)) & n<\infty \\
           f(x) & n=\infty
        \end{cases} \nonumber
    \end{align}

    \begin{claim}\label{claim:HLSClaim1}
        Suppose $f\in \ell^2(G_\infty)$ is pointwise non-negative.
        If $\supp(f)$ is compact, then $\supp(F)$ is compact.
    \end{claim}
\noindent \emph{Proof of \cref{claim:HLSClaim1}:} Since $\supp(f)\subseteq G_\infty$ is compact and $G_\infty$ is discrete, $\supp(f)=\set{z_1,\dots,z_m}$ for some $m\in \N$. Suppose for $n\in \N$, $F(n,q_n(x))>0$. We then have $\nnorm{f|_{xK_n}}_2>0$ and so $xK_n\cap \supp(f)\neq \varnothing$. From this we see $q_n(x)=q_n(z_i)$ for some $z_i\in \supp(f)$, which implies $(n,q_n(x))\in \ray_n(z_i)$. Therefore, $\supp(F)\subseteq \cup_{i=1}^m\ray_1(z_i)$, and $F$ is compactly supported by \Cref{lem:rays} Item \ref{lem-item:ray3}. 

\begin{claim}\label{claim:HLSClaim2}
    Suppose $f\in \ell^2(G_\infty)$ is pointwise non-negative.
    Then $F$ is continuous on $\G$, that is, for all $x\in G_\infty$,
    $$\lim_{n\to\infty} F(n,q_n(x))= f(x).$$
\end{claim}
\noindent \emph{Proof of \Cref{claim:HLSClaim2}:} Let $x\in G_\infty$ and $\epsilon>0$ be arbitrary. First note that
\begin{equation}\label{eq:HLSequation1}
F(n,q_n(x))^2=f_n(q_n(x))^2=\nnorm{f|_{xK_n}}_2^2=\sum_{z\in K_n}\abs{f(xz)}^2=\abs{f(x)}^2+\sum_{z\in K_n\setminus\set{e_{G_\infty}}}\abs{f(xz)}^2.    
\end{equation}
Since $\nnorm{f}_2^2<\infty$, there exists a finite set $E\subseteq G_\infty$ such that $\sum_{z\in G_\infty\setminus (E\cup\{x\})}\abs{f(z)}^2<\epsilon$.
This implies that, for every $x$,
\begin{align*}
    \sum_{y\in K_n\setminus(x^{-1}E\cup \set{e_{G_\infty}})}\abs{f(xy)}^2
    =\sum_{z\in xK_n\setminus(E\cup \set{x})}\abs{f(z)}^2
    \leq\sum_{z\in G_\infty\setminus (E\cup\{x\})}\abs{f(z)}^2<\epsilon. 
\end{align*}

Without loss of generality, assume $x\not\in E$.
        Let $x^{-1}E=\set{y_1,\dots,y_m}$. As $\set{K_n}_{n\in \N}$ is a nested sequence, $\cap_{n\in \N} K_n=\set{e_{G_\infty}}$, and $e_{G_\infty}\notin x^{-1}E$, there exists $N$ such that for all $n\geq N$, $y_i\notin K_n$ for $1\leq i\leq m$.
        So,  $x^{-1}E\cap K_n=\emptyset$ for $n\geq N$, and 
\begin{equation}\label{eq:HLSequation2}
    \sum_{z\in K_n\setminus \set{e_{G_\infty}}}\abs{f(xz)}^2=\sum_{z\in K_n\setminus (x^{-1}E\cup \set{e_{G_\infty}})}\abs{f(xz)}^2<\epsilon.
\end{equation}
        Putting \Cref{eq:HLSequation1} and \Cref{eq:HLSequation2} together, we have 
        $f(x)^2\leq F(n,q_n(x))^2\leq f(x)^2+\epsilon.$
        This implies that $\lim_{n\to \infty} F^2(n,q_n(x))= f^2(x)$.
        Since $f$ was assumed to be pointwise positive, we conclude that $F(n,q_n(x))\to F(\infty,x)=f(x)$ as $n\to \infty$.
        
    Note that given an arbitrary $f\in \ell^2(G_\infty)$, we may decompose $f$ into $f=f_1-f_2+i(f_3-f_4)$ where $f_i\geq 0$ pointwise for $1\leq i\leq 4$ by taking the positive and negative parts of the real and imaginary parts of $f$. For each $f_i$, construct $F_i$ as in \Cref{eq:HLSSection2} and define 
    \begin{align}\label{eq:HLSSection3}
        F:\G\to \C, \quad F=F_1-F_2+i(F_3-F_4). 
    \end{align}
    %
    %If we assume $f$ is continuous, then $f_i$ is continuous for each $i$. 
    Since each $F_i$ is continuous (by \Cref{claim:HLSClaim2}), the function $F$ is continuous. 
    Applying \Cref{claim:HLSClaim1}, if we take $f\in C_c(G_\infty)$ then $F\in C_c(\G)$ as well.

\begin{claim}\label{claim:HLSClaim3}
For every $f\in \ell^2(G_\infty)$,
the associated function $F$ is a continuous section of the left regular bundle $\Lambda$.
\end{claim}
\noindent \emph{Proof of \Cref{claim:HLSClaim3}:} 
Let $F$ be the function on $\G$ associated with $f\in \ell^2(G_\infty)$. The function $F\in C(\G)$ can be viewed as a section $\{F^n\}_{n\in {\mathbb N}\cup\{\infty\}}$, where $F^n=F|_{G_n}$, associated with the left regular bundle (\Cref{exa:L2(G)}).
For each $n\in {\mathbb N}$, let $S_n$ denote a set of representatives of quotients in $G_n$.

    First suppose that $f\geq 0$. Using the definition of continuous sections, we need to prove the following two items.
    \begin{enumerate}
        \item\label{Section:item1} Note that 
        $$\nnorm{F^n}_2^2=\sum_{x\in S_n}\abs{f_n(q_n(x))}^2=\sum_{x\in S_n}\nnorm{f|_{xK_n}}_2^2=\nnorm{f}_2^2.$$
        Constant functions are always continuous, and so the map $n\mapsto \nnorm{F^n}_2$ is continuous.
        \item\label{Section:item2} We show that $F$ is locally close to the fundamental sections of $\Lambda$, i.e., elements of $C_c(\G)$. Let $\epsilon>0$ be given, and choose $h\in C_c(G_\infty)$ with $\nnorm{h-f}_2<\epsilon$ and $h\geq 0$. Construct $H$ from $h$ via \Cref{eq:HLSSection2}. By \Cref{claim:HLSClaim1} and \Cref{claim:HLSClaim2}, $H\in C_c(\G)$. Note that
        $$\nnorm{F^n-H^n}_2^2=\sum_{x\in S_n}\abs{\nnorm{f|_{xK_n}}_2-\nnorm{h|_{xK_n}}_2}^2\leq \sum_{x\in S_n} \nnorm{f|_{xK_n}-h|_{xK_n}}_2^2=\nnorm{f-h}_2^2,$$
        which implies $\nnorm{F^n-H^n}_2<\epsilon$ for all $n\in \N$.
    \end{enumerate}
Next, consider a general $f\in \ell^2(G_\infty)$, and consider the decomposition $f=f_1-f_2+i(f_3-f_4)$, $f_i\geq 0$, and let $F=F_1-F_2+i(F_3-F_4)$ be defined as in \Cref{eq:HLSSection3}. By the previous case, each one of $F_1,F_2, F_3$ and $F_4$ are continuous sections of $\Lambda$. Since continuous sections form a complex linear space (\Cref{def:BanachField}), the linear combination $F=F_1-F_2+i(F_3-F_4)$ is a continuous section as well. This finishes the proof of the claim.
        
We have now demonstrated that for $f\in \ell^2(G_\infty)$, constructing $F$ as in \Cref{eq:HLSSection3} yields a continuous section of $\Lambda$ with $F^\infty=f$. Next, we show that this construction is norm-preserving in the following sense.
\begin{claim}\label{claim:HLSClaim4-norm}
For every $f\in \ell^2(G_\infty)$,
the associated continuous section $F$ satisfies $\|F\|_{\Delta}= \|f\|_2$.
\end{claim}
\noindent \emph{Proof of \Cref{claim:HLSClaim4-norm}:} 
Let $f=f_1-f_2+i(f_3-f_4)$ be the decomposition of $f$ into the positive and negative parts of the real and imaginary components of $f$, and $F=F_1-F_2+i(F_3-F_4)$ be defined as in \Cref{eq:HLSSection2}.  For $1\leq i\leq 4$, we have 
\begin{equation}\label{eq-norm-pos}
    \|F_i^n\|_{\ell^2(G_n)}^2=\sum_{x\in S_n} |F_i^n(q_n(x))|^2= \sum_{x\in S_n}\|f_i|_{xK_n}\|_2^2=\|f_i\|_2^2.
\end{equation}
Moreover, since by definition each function $F_i^n$ is pointwise non-negative,  we have
\begin{eqnarray*}
    \|F^n\|_{\ell^2(G_n)}^2&=&\|F_1^n-F_2^n+i(F_3^n-F_4^n)\|_{\ell^2(G_n)}^2=
    \|F_1^n-F_2^n\|_{\ell^2(G_n)}^2+\|F_3^n-F_4^n\|_{\ell^2(G_n)}^2\\
    &\leq& \|F_1^n\|_{\ell^2(G_n)}^2+\|F_2^n\|_{\ell^2(G_n)}^2+\|F_3^n\|_{\ell^2(G_n)}^2+\|F_4^n\|_{\ell^2(G_n)}^2.
\end{eqnarray*}
Using \Cref{eq-norm-pos} , we have
    
\begin{eqnarray*}
\|F^n\|_{\ell^2(G_n)}^2
&\leq& \|F_1^n\|_{\ell^2(G_n)}^2+\|F_2^n\|_{\ell^2(G_n)}^2+\|F_3^n\|_{\ell^2(G_n)}^2+\|F_4^n\|_{\ell^2(G_n)}^2\\
    &=& \|f_1\|_{\ell^2(G_\infty)}^2+\|f_2\|_{\ell^2(G_\infty)}^2+\|f_3\|_{\ell^2(G_\infty)}^2+\|f_4\|_{\ell^2(G_\infty)}^2\\
     &=& \|f_1-f_2\|_{\ell^2(G_\infty)}^2+\|f_3-f_4\|_{\ell^2(G_\infty)}^2\\
     &=&\|f\|_{\ell^2(G_\infty)}^2,
\end{eqnarray*}
where we used the fact that $f_1$ and $f_2$ (respectively $f_3$ and $f_4$) have disjoint support in the penultimate equality.
Since $F^\infty=f$, the section norm of $F$ can be computed as follows:
\begin{equation} \label{eq-section-norm}
    \|F\|_\Delta=\sup_{n\in \N\cup\{\infty\}} \|F^n\|_{\ell^2(G_n)}=\|f\|_{\ell^2(G_\infty)}.
\end{equation}
    
Let $u\in A(G_\infty)$ be arbitrary. Choose $f,g\in \ell^2(G_\infty)$ such that $u=L_{f,g}$ and $\nnorm{u}_{A(G_\infty)}=\nnorm{f}_2\nnorm{g}_2$ (\cite[Lemma 2.14]{Eymard}). Let $F,G$ be the sections constructed for $f,g$ as in \Cref{eq:HLSSection2} and \Cref{eq:HLSSection3}. It is clear that $\Lambda_{F,G}\in A(\G)$ and $R_\infty(\Lambda_{F,G})=L_{f,g}$. Since the restriction map
$R_\infty$ is contractive (\cite[Proposition 3.4 (iii)]{DeG}), we get
$$\nnorm{u}_{A(G_\infty)}=\nnorm{R_n(\Lambda_{F,G})}_{A(G_\infty)}\leq \nnorm{\Lambda_{F,G}}_{A(\G)}.$$
On the other hand, 
$$\nnorm{\Lambda_{F,G}}_{A(\G)}\leq \nnorm{F}_\Delta \nnorm{G}_\Delta =\nnorm{f}_2\nnorm{g}_2=\nnorm{u}_{A(G_\infty)}.$$
So $\nnorm{\Lambda_{F,G}}_{A(\G)}=\nnorm{u}_{A(G_\infty)}$, and the restriction map is a quotient map.
\end{proof}

\begin{rmk}
Here, we list two easy consequences of \Cref{thm:HLSSurjective}.
\begin{enumerate}
\item Positive-definite functions admit positive-definite norm-preserving lifts.
Indeed, suppose $u\in A(G_\infty)$ is positive-definite, and let $f\in\ell^2(G_\infty)$ be such that
$$ u(x)=\langle\lambda_{G_\infty}(x)f,f\rangle, \qquad \|f\|_2^2=u(e_{G_\infty})=\|u\|_{A(G_\infty)}. $$
Let $F$ be the continuous section lifted from $f$ as described in \Cref{eq:HLSSection2} and \Cref{eq:HLSSection3}. Then $v=\Lambda_{F,F}$ is a positive-definite element of $A(\G)$ satisfying
 $R_\infty(v)=u$ and $\|v\|_{A(\mathcal G)} = \|u\|_{A(G_\infty)}. $

\item 
The lifting property proved in \Cref{thm:HLSSurjective} can be extended to finitely many fibers.
More precisely, let $S\subseteq\mathbb N^\ast$ be a nonempty finite set. 
For each $n\in S$, let $u_n\in A(G_n)$ be given. 
Then there exists $v\in A(\mathcal G)$ such that 
\begin{equation}\label{eq:several-sec}
    R_n(v)=u_n \ \mbox{ and }\  \|v\|_{A(\G)}=\max_{n\in S}\|u_n\|_{A(G_n)}\ \forall n\in S.
\end{equation}
To see this, for each $n\in S$, let $f_n\in \ell^2(G_n)$ be chosen such that $u_n(\cdot)=\langle\lambda_{G_n}(\cdot) f_n,f_n\rangle$ with $\|u_n\|_{A(G_n)}=\|f_n\|_2^2$.
If $\infty\in S$, let $F=\{F_n\}_{n\in\N^*}$ denote the continuous section that lifts $f_\infty$ (described in \Cref{eq:HLSSection2} and \Cref{eq:HLSSection3}); otherwise, let $F$ denote the constant 0 section.
We now define the new continuous section $H=\{H_n\}_{n\in\N^*}$ by adjusting $F$ to fit the rest of $S$ as follows:
$$H_n=\left\{\begin{array}{ll}
  f_n   & n\in S\setminus\{\infty\} \\
 F_n    & \mbox{otherwise}
\end{array}\right..$$
It is easy to see that $v=\Lambda_{H,H}$ satisfies the conditions of \Cref{eq:several-sec}.
\end{enumerate}
\end{rmk}

\section{Applications: Fourier norm estimates and bounded approximate identity}
\label{sec:application}
The constructions in \Cref{thm:HLSSurjective} enable us to extend the result of \Cref{prop:normconvergence1}, Item \Cref{item:when-u-zero} to arbitrary elements of $A(\G)$. 
This extension yields a formula expressing the Fourier norm of $u\in A(\G)$ in terms of the norms of its restrictions to the fibers. We also characterize the HLS groupoids $\G$ for which $A(\G)$ admits a bounded approximate identity.

\begin{prop}\label{cor:norm-convergence-fibers-kernel}
    For every $u\in A(\G)$, we have $\lim_{n\to\infty}\nnorm{R_n(u)}_{A(G_n)}=\nnorm{R_\infty(u)}_{A(G_\infty)}$.
    Consequently, we can identify the kernel of the restriction map $R_\infty$ as follows:
    $$\ker(R_\infty)=\left\{u\in A(\G):  \|R_n(u)\|_{A(G_n)}\to 0\right\}.$$
We denote the subspace defined on the right hand side by  $c_0-\bigoplus_{n\in\mathbb N}A(G_n)$ .
\end{prop}
\begin{proof}
First, fix $u\in A(\G)\cap C_c(\G)$. Applying \Cref{thm:HLSSurjective} to $R_\infty(u)$, there exists $v\in A(\G)$ such that
\[
R_\infty(v)=R_\infty(u).
\]
Moreover, \Cref{claim:HLSClaim1} implies that $v\in C_c(\G)$. Thus, $u-v\in C_c(\G)$ and
$R_\infty(u-v)=0.$
Applying \Cref{prop:normconvergence1} \eqref{item:when-u-zero} to $u-v$, we obtain
$\lim_{n\to\infty}\nnorm{R_n(u)-R_n(v)}_{A(G_n)}=0.$
Consequently,
\begin{equation}\label{eq:vn-un}
    \lim_{n\to\infty}
\abs{\nnorm{R_n(u)}_{A(G_n)}-\nnorm{R_n(v)}_{A(G_n)}}=0.
\end{equation}

Next, we show that
$\lim_{n\to\infty}\nnorm{R_n(v)}_{A(G_n)}
=
\nnorm{R_\infty(v)}_{A(G_\infty)}.$
First, note that $\nnorm{v}_{A(\G)}=\nnorm{R_\infty(v)}_{A(G_\infty)}$ (\Cref{thm:HLSSurjective}). 
So, by \Cref{prop:normconvergence1} \eqref{item:lower-bnd-gen},
$\limsup_{n\to\infty}\nnorm{R_n(v)}_{A(G_n)}
\leq
\nnorm{R_\infty(v)}_{A(G_\infty)}.$
On the other hand, \Cref{prop:normconvergence1} \eqref{item:upper-bnd-A-cf} gives
\[
\nnorm{R_\infty(v)}_{A(G_\infty)}
\leq
\liminf_{n\to\infty}\nnorm{R_n(v)}_{A(G_n)}.
\]
It follows that
$\limsup_{n\to\infty}\nnorm{R_n(v)}_{A(G_n)}
\leq
\nnorm{R_\infty(v)}_{A(G_\infty)}
\leq
\liminf_{n\to\infty}\nnorm{R_n(v)}_{A(G_n)},
$
and therefore
\begin{equation}\label{eq:norm-v}
\lim_{n\to\infty}\nnorm{R_n(v)}_{A(G_n)} = \nnorm{R_\infty(v)}_{A(G_\infty)}.
\end{equation}
Combining \Cref{eq:vn-un} and \Cref{eq:norm-v}, together with the fact that $R_\infty(v)=R_\infty(u)$, 
we get that
\[
\lim_{n\to\infty}\nnorm{R_n(u)}_{A(G_n)}
=
\nnorm{R_\infty(u)}_{A(G_\infty)}
\]
for every $u\in A(\G)\cap C_c(\G)$.

Finally, let $u\in A(\G)$ and fix $\epsilon>0$. Choose
$v\in A(\G)\cap C_c(\G)$ such that $\nnorm{u-v}_{A(\G)}<\frac{\epsilon}{3}.$
Since each restriction map is contractive, for every $n\in\N^*$ we have
\[
\abs{\nnorm{R_n(u)}_{A(G_n)}-\nnorm{R_n(v)}_{A(G_n)}}
\leq
\nnorm{R_n(u-v)}_{A(G_n)}
\leq
\nnorm{u-v}_{A(\G)}
<
\frac{\epsilon}{3}.
\]
By the first part of the proof, there exists $N\in\N$ such that, for every $n\geq N$,
\[
\abs{\nnorm{R_n(v)}_{A(G_n)}
-\nnorm{R_\infty(v)}_{A(G_\infty)}}
<
\frac{\epsilon}{3}.
\]
Therefore, for every $n\geq N$,
\begin{align*}
\abs{\nnorm{R_n(u)}_{A(G_n)}
-\nnorm{R_\infty(u)}_{A(G_\infty)}}
&\leq
\abs{\nnorm{R_n(u)}_{A(G_n)}
-\nnorm{R_n(v)}_{A(G_n)}}\\
&\quad+
\abs{\nnorm{R_n(v)}_{A(G_n)}
-\nnorm{R_\infty(v)}_{A(G_\infty)}}\\
&\quad+
\abs{\nnorm{R_\infty(v)}_{A(G_\infty)}
-\nnorm{R_\infty(u)}_{A(G_\infty)}}\\
&<\epsilon.
\end{align*}
Hence,
$\nnorm{R_n(u)}_{A(G_n)}
\longrightarrow
\nnorm{R_\infty(u)}_{A(G_\infty)}.$

To prove the second statement, first observe that 
$c_0\text{-}\bigoplus_{n\in\mathbb N} A(G_n)\subseteq \ker(R_\infty).$
For the reverse inclusion, let $u\in \ker(R_\infty)$. By the first statement of the proposition,
$$\lim_{n\to\infty}\|R_n(u)\|_{A(G_n)}=\|R_\infty(u)\|_{A(G_\infty)}=0.$$
Thus, the sequence $\{R_n(u)\}_{n\in\mathbb N}$ belongs to
$c_0\text{-}\bigoplus_{n\in\mathbb N} A(G_n),$
and we are done.
\end{proof}

\begin{cor}\label{cor:norm-A(G)-limit-sections}
    For every $u\in A(\G)$, we have $\nnorm{u}_{A(\G)}=\sup_{n\in \N^*}\nnorm{R_n(u)}_{A(G_n)}$.
\end{cor}
\begin{proof}
    By \Cref{prop:normconvergence1} \eqref{item:lower-bnd-gen}, $\nnorm{u}_{A(\G)}\geq \sup\limits_{n\in \infN}\nnorm{R_n(u)}_{A(G_n)}.$ So we only need to show the reverse inequality.
    
    First assume that $u\in A(\G)$ satisfies $R_\infty(u)=0$. 
    By \Cref{cor:norm-convergence-fibers-kernel}, $\lim_{n\to\infty}\|R_n(u)\|_{A(G_n)}=0$. 
    For each $n\in\N$ we may choose $f_n,h_n\in\ell^2(G_n)$ such that
    \[R_n(u)=\langle \lambda_{G_n}(\cdot)f_n,h_n\rangle_{\ell^2(G_n)} \quad\text{and}\quad \|f_n\|_2=\|h_n\|_2=\|R_n(u)\|_{A(G_n)}^{1/2}.\]
    Set $f_\infty=g_\infty=0$, and let $F=\{f_n\}_{n\in\N^*}$ and $H=\{h_n\}_{n\in\N^*}$. 
    By \Cref{cor-cts-sec} \eqref{item:section2} both $F$ and $H$ belong to $\Gamma$, 
    and by \Cref{cor-cts-sec} \eqref{item:sections-coefficients}, we have $\Lambda_{F,H}\in A(\mathcal{G})$.
    Since $R_n(\Lambda_{F,H})=R_n(u)$ for every $n\in \N^*$, we have $\Lambda_{F,H}=u$. Therefore,
\[
  \|u\|_{A(\mathcal{G})}
  \leq\|F\|_\Delta\,\|H\|_\Delta
  =\Big(\sup_{n\in\N}\|f_n\|_2\Big)\Big(\sup_{n\in\N}\|h_n\|_2\Big)
  =\sup_{n\in\mathbb{N}}\|R_n(u)\|_{A(G_n)}.
\]
This proves the desired statement for the case where $R_\infty(u)=0$.

Next, let $u\in A(\G)$ be arbitrary and put
$S:=\sup_{n\in\N^*}\|R_n(u)\|_{A(G_n)}$, which is finite by
\Cref{prop:normconvergence1} \eqref{item:lower-bnd-gen}.  
For each
$n\in\N^*$ choose $\xi_n,\eta_n\in\ell^2(G_n)$ such that
\begin{equation*}\label{eq:xi-eta}
  R_n(u)=\langle\lambda_{G_n}(\cdot)\xi_n,\eta_n\rangle
  \quad\text{and}\quad
  \|\xi_n\|_2^2=\|\eta_n\|_2^2=\|R_n(u)\|_{A(G_n)}\le S .
\end{equation*}
Let $F=\{f_n\}_{n\in\N^*},H=\{h_n\}_{n\in \N^*}$ be the continuous sections associated with
$\xi_\infty$ and $\eta_\infty$ defined by \Cref{eq:HLSSection2} and \Cref{eq:HLSSection3} (\Cref{claim:HLSClaim2}). Clearly,
$f_\infty=\xi_\infty$, $h_\infty=\eta_\infty$, and by
\Cref{eq-section-norm},
\begin{equation}\label{eq:lift-bound}
  \|f_n\|_2\le\|\xi_\infty\|_2
  \quad\text{and}\quad
  \|h_n\|_2\le\|\eta_\infty\|_2
  \qquad(n\in\mathbb{N}^*).
\end{equation}
Let $v:=\Lambda_{F,H}\in A(\mathcal{G})$. 
Clearly, $R_\infty(u-v)=0$, and by \Cref{cor:norm-convergence-fibers-kernel},
$\|R_n(u-v)\|_{A(G_n)}\to 0$.
Let $\epsilon>0$ be arbitrary, and fix $N\in\N$ such that
\begin{equation}\label{eq:choice-N}
  \|R_n(u-v)\|_{A(G_n)}<\epsilon\qquad\text{for all } n\ge N .
\end{equation}
Next, we adjust the continuous sections $F$ and $H$ so that the corresponding coefficient function matches $u$ on the first $N-1$ fibers.
Namely, we define $F'=\{f'_n\}_{n\in \N^*}$ and $H'=\{h'_n\}_{n\in \N^*}$ by
\[
  f'_n:=
  \begin{cases}
    \xi_n, & n<N,\\
    f_n, & n\in[N,\infty],
  \end{cases}
  \qquad
  h'_n:=
  \begin{cases}
    \eta_n, & n<N,\\
  h_n, & n\in[N,\infty].
  \end{cases}
\]
Clearly, $F'$ and $H'$ are continuous sections, as they eventually match $F$ and $H$ respectively.
Let $w:=\Lambda_{F',H'}$. By \Cref{cor-cts-sec} \eqref{item:sections-coefficients}
$w\in A(\mathcal{G})$, and by \Cref{eq:lift-bound}
\begin{equation}\label{eq:w-bound}
  \|w\|_{A(\mathcal{G})}\le\|F'\|_\Delta\,\|H'\|_\Delta\le S .
\end{equation}
Next, we estimate $\|u-w\|_{A(\G)}$. By construction of $w$, we have $R_\infty(u-w)=0$. So, by the first part of the proof, we have 
\begin{align*}
   \|u-w\|_{A(\G)}&=\sup_{n\in \N^*}\|R_n(u-w)\|_{A(G_n)}\\
   =&\sup_{n\in[N,\infty]}\|R_n(u)-R_n(w)\|_{A(G_n)}  \\
   =&\sup_{n\in[N,\infty]}\|R_n(u)-R_n(v)\|_{A(G_n)}, 
\end{align*}
so by \Cref{eq:choice-N}, we have $\|u-w\|_{A(\G)}\leq \epsilon$. Putting this fact together with \Cref{eq:w-bound} completes the proof, as for each $\epsilon>0$, we have
$$\|u\|_{A(\G)}\leq \|u-w\|_{A(\G)}+\|w\|_{A(\G)}\leq S+\epsilon.$$
\end{proof}

\begin{rmk}\label{rem:norm-A-obtained-Lambda}
For a locally compact group, the norm of an element of $A(G)$ can be obtained using only coefficient functions of the left regular representation. An analogous result for groupoids would, in particular, allow the Fourier algebra to be defined without first introducing the Fourier--Stieltjes algebra.
We were not able to extend this group property to groupoids (or even HLS groupoids) in general. 
However, as an easy corollary of the previous result, we observe that for every $u\in A(\G)$ with $R_\infty(u)=0$, we have
$$\nnorm{u}_{A(\G)}=\inf\set{\nnorm{F}_{\Delta}\nnorm{G}_{\Delta}:u=\Lambda_{F,G}}.$$
Indeed, in the proof of \Cref{cor:norm-A(G)-limit-sections}, we constructed sections $F,H\in \Gamma$ such that $\Lambda_{F,H}=u$, and 
\[
  \|u\|_{A(\mathcal{G})}
  \leq\|F\|_\Delta\,\|H\|_\Delta
  =\sup_{n\in\mathbb{N}}\|R_n(u)\|_{A(G_n)}\leq   \|u\|_{A(\mathcal{G})},\]
where the final inequality follows from contractility of restriction maps $R_n$. 
Thus, equality holds throughout, proving the desired statement.
\end{rmk}

We finish this section by presenting another interesting applications of \Cref{thm:HLSSurjective} and \Cref{cor:norm-A(G)-limit-sections} as follows.
We first recall one more definition from Banach algebra theory.
\begin{defn}\label{def:bai}
Let $\mathcal{A}$ be a Banach algebra. A \emph{left} (respectively \emph{right})
\emph{approximate identity} for $\mathcal{A}$ is a net $\{e_\alpha\}_{\alpha\in I}$
in $\mathcal{A}$ such that
\[
  \lim_{\alpha}\|e_\alpha a - a\| = 0
  \qquad\Big(\text{respectively } \lim_{\alpha}\|a e_\alpha - a\| = 0\Big)
  \qquad\text{for all } a\in\mathcal{A}.
\]
A net that is both a left and a right approximate identity is called a
(\emph{two-sided}) \emph{approximate identity}. An approximate identity
$\{e_\alpha\}_{\alpha\in I}$ is \emph{bounded} if there is a constant $C\geq 1$ with
\[
  \sup_{\alpha\in I}\|e_\alpha\| \leq C .
\]
In this case $\{e_\alpha\}_{\alpha\in I}$ is called a \emph{bounded approximate
identity} with bound $C$. 
\end{defn}
\begin{cor}[Leptin's theorem for HLS groupoids]
\label{cor:BAI}
    Let $\G:=\sqcup_{n\in {\mathbb N}^*} G_n$ be an HLS groupoid. Then $A(\G)$ has a bounded approximate identity if and only if $G_\infty$ is amenable. Moreover, in this case $A(\G)$ has a bounded approximate identity bounded by $1$.
\end{cor}
\begin{proof}
     Suppose that $\{e_\alpha\}_{\alpha\in A}$ is a bounded approximate identity for $A(\G)$. Since $R_\infty$ is a  contractive algebra homomorphism, $\{R_\infty(e_\alpha)\}_{\alpha\in A}$ is a bounded approximate identity for $A(G_\infty)$, and $G_\infty$ is amenable by Leptin's theorem \cite{Leptin}. This proves the forward direction.

     To prove the backward direction,  suppose that $G_\infty$ is amenable. By Leptin's theorem \cite{Leptin}, $A(G_\infty)$ has a bounded approximate identity $\{\varphi_\beta\}_{\beta\in B}$ with $\nnorm{\varphi_\beta}_{A(G_\infty)}\leq 1$ for all $\beta\in B$. 
     By \Cref{thm:HLSSurjective}, for each $\beta\in B$ there exists $v_\beta\in A(\G)$ such that
    \[
        R_\infty(v_\beta)=\varphi_\beta
        \quad\text{and}\quad
        \nnorm{v_\beta}_{A(\G)}=\nnorm{\varphi_\beta}_{A(G_\infty)}\leq 1.
    \]
    Next, we define a new net of functions in $A(\G)$ indexed over the directed set $B\times \N$ that will act as a bounded approximate identity for $A(\G)$.
    For $\beta\in B$ and $N\in\N$, define
    \[
        e_{\beta,N}:=v_\beta+\sum_{n=1}^{N-1}\iota_n\big(1_{G_n}-R_n(v_\beta)\big),
    \]
    where $\iota_n:A(G_n)\to A(\G)$ is the embedding of \Cref{lem:spectrumembedding}, and $1_{G_n}$ denotes the constant function on $G_n$ attaining 1 everywhere. Since $G_n$ is finite, we know $1_{G_n}\in A(G_n)$ and $\|1_{G_n}\|_{A(G_n)}=1$. 
    So, $e_{\beta,N}\in A(\G)$. Moreover, applying \Cref{cor:norm-A(G)-limit-sections}, we get
    $$\|e_{\beta,N}\|_{A(\G)}=\max\left\{\max_{n\in [1, N-1]}\|1_{G_n}\|_{A(G_n)},\sup_{n\in [N,\infty]}\|R_n(v_\beta)\|_{A(G_n)}\right\}\leq 1,$$
    since $\|v_\beta\|_{A(\G)}\leq 1$.
    So $\{e_{\beta,N}\}_{(\beta,N)\in B\times \N}$ is a bounded net of functions in $A(\G)$, with bound 1.
    
    To complete the proof, it only remains to show that this net is an approximate identity for $A(\G)$.
    Let $u\in A(\G)$ be arbitrary. First note that by \Cref{cor:norm-A(G)-limit-sections} and the definition of $e_{\beta,N}$, for every $\beta\in B$, and $N\in\N$ we have
    \begin{equation}\label{eq:norm-of-difference}
        \nnorm{u\,e_{\beta,N}-u}_{A(\G)}=\sup_{n\in[N,\infty]}\nnorm{R_n(u \, e_{\beta,N}-u)}_{A(G_n)}.
    \end{equation}
    Let $\epsilon>0$ be arbitrary.
    Since $\{\varphi_\beta\}_{\beta\in B}$ is an approximate identity for $A(G_\infty)$, there exists $\beta_0\in B$ such that for all $\beta\geq \beta_0$ and every $N\in \N$, we have
    \begin{equation}\label{eq:last-fiber-approx}
        \nnorm{R_\infty(u\, e_{\beta, N}-u)}_{A(G_\infty)}=\nnorm{R_\infty(u)\varphi_\beta-R_\infty(u)}_{A(G_\infty)}<\frac{\epsilon}{2},
    \end{equation}
    where in the first equality we used the fact that $\varphi_\beta=R_\infty(v_\beta)=R_\infty(e_{\beta, N})$.
    Applying \Cref{cor:norm-convergence-fibers-kernel} to $u\, e_{\beta, 1}-u$, there exists $N_0\in \N$ such that for all $n\geq N_0$, we have
     \begin{equation}\label{eq:limit-of-fibers}
        \nnorm{R_n(u\, e_{\beta, 1}-u)}_{A(G_n)}<\epsilon.
    \end{equation}
    Combining \Cref{eq:norm-of-difference}, \Cref{eq:last-fiber-approx} and \Cref{eq:limit-of-fibers}, we get
    $\|u\, e_{\beta, N}-u\|_{A(\G)}<\epsilon$ for all $N\geq N_0$ and $\beta\geq \beta_0$, so $\{u_{\beta,n}\}_{(\beta,n)\in B\times \N}$ forms an approximate identity for $A(\G)$.
\end{proof}

%%%%%%%%%%%%%%%%%%%%%%%%%%%%%%%%%%%%%%%%%%%%%%%%%%%%%%%%%%%%%%%%%%%%%%%%%%%%%%%%%%%%%%%%%%%%%%%%%%%%%%%%%%%%%%%%%%%%%%%%%%%%%%%%%%%%%%%%%%%%%%%%%%%%%%%%%%%%%%%%%%%%%%%%%%%%%%%%%%%%%%%%%%%%%%%%%%%%%%%%%%%%%%%%%%%%%%%%%%%%%%%%%%%%%%%%%%%%%%
\section{Gelfand spectrum of the Fourier algebra}
\label{sec:gelfandspectrum}
Let $G$ be a locally compact group. The Gelfand spectrum of $A(G)$, denoted $\sigma(A(G))$, is the set of all nonzero multiplicative linear functionals on $A(G)$. A classical result of Eymard shows that
\[
\sigma(A(G)) \cong G
\]
homeomorphically \cite[Theorem 3.34]{Eymard}. In \cite{PatersonPaper}, Paterson remarks that it is unknown whether the analogous statement holds for arbitrary groupoids. 
The discrete groupoid case was established in \cite{dissertation}. 
In this section, we answer Paterson's question for HLS groupoids.

\begin{thm}\label{prop:spectrum}
     Let $\G:=\sqcup_{n\in {\mathbb N}^*} G_n$ be an HLS groupoid. Then $\sigma(A(\G))$ is homeomorphic to $\G$.
\end{thm}
\begin{proof}
Note that $\G$ is a locally compact and Hausdorff space, and  $A(\G)$ is a subalgebra of $C_0(\G)$ that strongly separates points of $\G$, i.e., for every distinct $x,y\in \G$ there exists $f\in A(\G)$ with $f(x)\neq f(y)$. 
Thus, by \cite[Theorem 3.2.4]{Rickart}, the natural embedding $\G\to \sigma(A(\G))$, given by point evaluation, 
is a homeomorphism into $\sigma(A(\G))$. What remains is to demonstrate that any multiplicative, linear functional on $A(\G)$ must correspond to point evaluation for a given groupoid element.

To that end, let $\Phi:A(\G)\to \C$ be a non-zero, bounded, multiplicative linear functional. For each $n\in \N$, let $\iota_n: A(G_n)\hookrightarrow A(\G)$ be the embedding as in \Cref{lem:spectrumembedding}.  
First, we observe that $\Phi|_{\iota_n(A(G_n))}\neq 0$ for at most one $n\in \N$. 
Indeed, if there were distinct $n,m$ such that $\Phi(u)\neq 0$ and $\Phi(v)\neq 0$ for $u\in \iota_n(A(G_n))$ and $v\in \iota_m(A(G_m))$, then since $uv=0$, we necessarily have
$$0=\Phi(0)=\Phi(uv)=\Phi(u)\Phi(v)$$
a contradiction.

To prove the result, we consider two cases:

\emph{Case 1:} There exists exactly one $n_0\in \N$ such that for every $n\in \N$, $n\neq n_0$, $\Phi|_{\iota_n(A(G_n))}=0$, and $\Phi(\iota_{n_0}(v))\neq 0$ for some $v\in A(G_{n_0})$. Define $\psi:A(G_{n_0})\to \C$ by $\psi(u)=\Phi(\iota_{n_0}(u))$. Clearly, $\psi$ is a nonzero bounded, multiplicative, linear functional on $A(G_{n_0})$. By \cite[Theorem 3.34]{Eymard}, there exists some $x_0\in G_{n_0}$ such that $\psi(w)=w(x_0)$ for every $w\in A(G_{n_0})$.

\begin{claim}\label{claim:Case1HLSSPectrum}
    If $u\in A(\G)$ with $u|_{G_{n_0}}=0$, then $\Phi(u)=0$.
\end{claim}
\noindent \emph{Proof of \Cref{claim:Case1HLSSPectrum}:} Recall that $v\in A(G_{n_0})$ was chosen such that $\Phi(\iota_{n_0}(v))\neq 0$. 
Clearly $u\, \iota_{n_0}(v)=0$, so we have 
$$0=\Phi(0)=\Phi(\iota_{n_0}(v)u)=\Phi(\iota_{n_0}(v))\Phi(u).$$ Since $\Phi(\iota_{n_0}(v))$ is non-zero, we must have that $\Phi(u)=0$, demonstrating the claim.

To conclude this case, we decompose $f\in A(\G)$ in the form
$f=f-\iota_{{n_0}}(R_{n_0}(f))+\iota_{{n_0}}(R_{n_0}(f))$. 
So, 
$$\Phi(f)=\Phi(f-\iota_{{n_0}}(R_{n_0}(f)))+\Phi(\iota_{{n_0}}(R_{n_0}(f)))=\Phi(f-\iota_{{n_0}}(R_{n_0}(f)))+\psi(R_{n_0}(f)).$$ 
The first term is 0 by \Cref{claim:Case1HLSSPectrum} and the second one is $f(n_0,x_0)$, proving that $\Phi$ is just the evaluation functional at the point $(n_0,x_0)\in \G$. This proves the result for Case 1.
\medskip

\emph{Case 2:} Suppose that for every $n\in \N$, $\Phi|_{\iota_{n}(A(G_n))}=0$. 
Let $\mathds{1}$ denote the constant function 1 on $\G$, and $\mathds{1}_n:\G\to \C$ be the characteristic function of $G_n$.
Let $u\in A(\G)$ be an arbitrary element. 
It is  easy to see that $\iota_k(R_k(u))=\mathds{1}_ku$.
Using \Cref{lem:spectrumembedding} and the fact that each $R_k:A(\G)\to A(G_k)$ is a homomorphism, for every $n\geq 1$, the functions $\mathds{1}_n$ and  $u_n=\sum_{k=1}^n\iota_k(R_k(u))$ belong to $A(\G)$. 
Since $\Phi|_{\iota_{k}(A(G_k))}=0$ for $k\in \N$, we have 
\begin{equation}\label{eq:Phi:u-un}
    \Phi(u)=\Phi(u-u_n)+\Phi\of{\sum_{k=1}^n\iota_k(R_k(u))}=\Phi(u-u_n), \ \forall n\in \N.
\end{equation}

Since $A(\G)\cap C_c(\G)$ is dense in $A(\G)$ and $\Phi$ is nonzero, there exists a compactly supported $v\in A(\G)$ satisfying $\Phi(v)\neq 0$. By \Cref{lem:rays}, Items \Cref{lem-item:ray3} and \Cref{lem-item:ray2}, there exist $x_1,\ldots,x_m\in G_\infty$ and $N\in \N$ such that
\begin{equation*}
    \supp(v-v_N)\subseteq \sqcup_{i=1}^m\ray_N(x_i),
\end{equation*}
where $v_N=\sum_{n=1}^N\mathds{1}_nv$.
Using \eqref{eq:Phi:u-un}, we have $\Phi(v-v_N)\neq 0.$
So, replacing $v$ by $v-v_N$ if necessary, there exists $v\in A(\G)$ satisfying the following:
\begin{equation}
    v\in A(\G)\cap C_c(\G), \quad \Phi(v)\neq 0, \quad \supp(v)\subseteq \sqcup_{i=1}^m\ray_N(x_i).
\end{equation}
Moreover, $v=v{\mathds 1}_{\ray_N(x_1)}+\ldots+v{\mathds 1}_{\ray_N(x_m)}$, and each summand belongs to $A(\G)$ by  \Cref{cor:characteristic-in-A}. Since $\Phi(v{\mathds 1}_{\ray_N(x_1)})+\ldots+\Phi(v{\mathds 1}_{\ray_N(x_m)})\neq 0$, without loss of generality we can assume that 
$\Phi(v{\mathds 1}_{\ray_N(x_1)})\neq 0$.
Again replacing $v$ by $v{\mathds 1}_{\ray_N(x_1)}$ if necessary, we can assume that $v$ is supported in the single ray $\ray_N(x_1)$.

\begin{claim}\label{claim:phi-of-rays}
$\Phi({\mathds 1}_{\ray_1(x_1)})=1$ and $\Phi({\mathds 1}_{\ray_1(x)})=0$ for $x\neq x_1$.
Consequently, if $u\in A(\G)\cap C_c(\G)$ and $(\infty, x_1)\not\in \supp(u)$ then $\Phi(u)=0$.
\end{claim}
\noindent \emph{Proof of \Cref{claim:phi-of-rays}:} 
Since $v$ is supported in  $\ray_1(x_1)$, we have that 
$$0\neq \Phi(v)=\Phi(v{\mathds 1}_{\ray_1(x_1)})=\Phi(v)\Phi({\mathds 1}_{\ray_1(x_1)}),$$
which implies that $\Phi({\mathds 1}_{\ray_1(x_1)})=1$. Moreover, for every $x\in G_\infty$ with $x\neq x_1$, there exits $K\in \N$ for which we have  ${\ray_K(x)}\cap {\ray_K(x_1)}=\emptyset$, implying that $(v-v_K){\mathds 1}_{\ray_K(x_1)}=0$.
So, 
$$0=\Phi((v-v_K){\mathds 1}_{\ray_K(x_1)})=\Phi(v-v_K)\Phi({\mathds 1}_{\ray_K(x)})=\Phi(v)\Phi({\mathds 1}_{\ray_1(x)}),$$
where we used \eqref{eq:Phi:u-un} in the last equality. Since $\Phi(v)\neq 0$, we must have 
$\Phi({\mathds 1}_{\ray_1(x)}).$

Next, let $u\in A(\G)\cap C_c(\G)$ be such that $(\infty,x_1)\not\in \supp(u)$. Applying \Cref{lem:rays} again,
there exists $N'\in \N$ such that $u-u_{N'}$ is supported in finitely many disjoint rays whose intersection with the infinite fiber $G_\infty$ lies inside $\supp(u)$. 
That is, there exist $z_1,\ldots, z_t\in G_\infty$ so that
$$\supp(u-u_{N'})\subseteq \sqcup_{i=1}^t\ray_{N'}(z_i)\ \mbox{ and } \ (\infty,z_i)\in \supp(u).$$
Consequently $u-u_{N'}= (u-u_{N'})({\mathds 1}_{\ray_1(z_1)}+\ldots+{\mathds 1}_{\ray_1(z_t)})$ and $z_i\neq x_1$. So, using the first part of the claim, we have
$$\Phi(u)=\Phi(u_{N'})+\Phi(u-u_N')\left(\Phi({\mathds 1}_{\ray_1(z_1)})+\ldots+\Phi({\mathds 1}_{\ray_1(z_t)})\right)=0,$$
which finishes the proof of the claim.

\begin{claim}\label{claim:u-supports-x1}
If $u\in A(\G)\cap C_c(\G)$ and $(\infty,x_1)\in \supp(u)$ then $\Phi(u)=u(\infty,x_1)$.
\end{claim}
\noindent \emph{Proof of \Cref{claim:u-supports-x1}:} 
Let $\widetilde{u}=u-u(\infty,x_1){\mathds 1}_{\ray_1(x_1)}$.
Clearly, $\widetilde{u}$ is a compactly supported element of $A(\G)$ as well. Assuming $(\infty,x_1)\notin\supp(\widetilde{u})$, by \Cref{claim:phi-of-rays} we have $\Phi(\widetilde{u})=0$. This finishes the proof of the current claim under this assumption, since 
$$\Phi(u)=u(\infty,x_1)\Phi({\mathds 1}_{\ray_1(x_1)})=u(\infty,x_1).$$
Next, suppose $(\infty,x_1)\in\supp(\widetilde{u})$.
By the definition of $\widetilde{u}$, we have $\widetilde{u}(\infty, x_1)=0$. 
Since $\widetilde{u}$ is compactly supported, there exist $M\in \N$ and $y_1,\ldots,y_s\in G_\infty$ satisfying \begin{equation}
    y_1=x_1\quad \mbox{ and }\quad \supp(\widetilde{u}-\widetilde{u}_M)\subseteq \sqcup_{i=1}^s{\mathds 1}_{\ray_M(y_i)}.
\end{equation}
Therefore, $\widetilde{u}-\widetilde{u}_M=(\widetilde{u}-\widetilde{u}_M)(\sum_{i=1}^s{\mathds 1}_{\ray_M(y_i)})$, and applying Equation \eqref{eq:Phi:u-un} and \Cref{claim:phi-of-rays} again we have
\begin{equation}\label{eq:eq:u-tlde-cut}
    \Phi(\widetilde{u})=\Phi(\widetilde{u}-\widetilde{u}_M)=\Phi((\widetilde{u}-\widetilde{u}_M){\mathds 1}_{\ray_M(x_1)}).
\end{equation}
For each $m\geq M$, define $w_m:=(\widetilde{u}-\widetilde{u}_m){\mathds 1}_{\ray_m(x_1)}$, and note that $w_m\in A(\G)$.
Since $w_m$ is supported in $\ray_m(x_1)$, we can write it as a coefficient function as follows:
$$w_m=\Lambda_{\bxi, \boldeta}, \ \mbox{ where } 
\xi_n:=\left\{\begin{array}{cc}
   (u(n,q_n(x_1))-u(\infty,x_1))\delta_{q_n(x_1)}  & n\geq m \\
    0 & n<m
\end{array}\right.,
\ \eta_n:= \left\{\begin{array}{cc}
\delta_{q_n(e_\infty)} & n\geq m \\
    0 & n<m
\end{array}\right..$$
Clearly $\boldeta\in \Delta$. The fact that $\bxi\in \Delta$ follows from the continuity of $u$ at $(\infty,x_1)$.
From the above representation of $w_m$, we get
$$\|w_m\|_{A(\G)}\leq \|\bxi\|_\Delta\|\boldeta\|_\Delta=\sup_{n\geq m}|u(n,q_n(x_1))-u(\infty,x_1)|.$$
Combining \eqref{eq:Phi:u-un} and \eqref{eq:eq:u-tlde-cut}, we get
$$|\Phi(\widetilde{u})|=\lim_{m\to \infty}|\Phi(w_m)|\leq \lim_{n\to \infty}|u(n,q_n(x_1))-u(\infty,x_1)|=0.$$
This completes the proof of the claim, as by \Cref{claim:phi-of-rays}, we have $\Phi({\mathds 1}_{\ray_1(x_1)})=1$.

Putting \Cref{claim:phi-of-rays} and \Cref{claim:u-supports-x1} together concludes Case 2. So far we have proved that there exist $n_1\in \N\cup\{\infty\}$ and $x_1\in G_{n_1}$ satisfying the following:
\begin{eqnarray*}
    \forall u\in A(\G)\cap C_c(\G),\quad \Phi(u)=u(n_1,q_{n_1}(x_1)). 
\end{eqnarray*}
It is easy to see that the same equality holds for every element of $A(\G)$. Indeed, let $u\in A(\G)$ and consider a sequence $\{f_n\}_n\subseteq A(\G)\cap C_c(\G)$ with $\|u-f_n\|_{A(\G)}\to 0$ as $n\to \infty$.
By continuity of $\Phi$, we have
$$\Phi(u)=\lim_{n\to \infty}\Phi(f_n)=\lim_{n\to \infty} f_n(n_1,q_{n_1}(x_1))=u(n_1,q_{n_1}(x_1)),$$
where the last equality holds because convergence in $A(\G)$ implies pointwise convergence.
\end{proof}

%%%%%%%%%%%%%%%%%%%%%%%%%%%%%%%%%%%%%%%%%%%%%%%%%%%%%%%%%%%%%%%%%%%%%%%%%%%%%%%%%%%%%%%%%%%%%%%%%%%%%%%%%%%%%%%%%%%%%%%%%%%%%%%%%%%%%%%%%%%%%%
\section{Extending beyond HLS groupoids}\label{section:beyond-HLS}
Some of the results established in this paper extend beyond HLS groupoids to arbitrary \'{e}tale groupoids. 
In this short section, we present these extensions. First, we prove an analogue of \Cref{thm:HLSSurjective} for all \'{e}tale groupoids. Unlike the proof of \Cref{thm:HLSSurjective}, however, the argument given here is non-constructive. 
Second, we show that the spectrum of $A(\G)$ is homeomorphic to $\G$ under the map given by point evaluation.
This generalizes the result of \Cref{prop:spectrum} and answers a question raised by Paterson in \cite[Section~6]{PatersonPaper} for the case of \'{e}tale groupoids.
However, Paterson's question remains open in full generality.

\begin{prop}\label{prop:etale-lift}
    Let $\G$ be a (locally compact Hausdorff) \'{e}tale groupoid, and consider $a\in\G^{(0)}$. 
    Let $R_a$ denote the restriction map $A(\G)\to A(\G^a_a)$.
    Then for each $u\in A(\G_a^a)$, there exists $v\in A(\mathcal G)$ such that $R_a(v)=u$ and $\|v\|_{A(\mathcal G)}=\|u\|_{A(\G_a^a)}$. 
\end{prop}
\begin{proof}
Consider a representation of $u$ as a coefficient function $\langle\lambda_G(\cdot)f,g\rangle$, where $f,g\in \ell^2(\G_a^a)$ satisfying $\|u\|_{A(\G_a^a)}=\|f\|_2\|g\|_2$.
Since $\G_a^a$ is a subset of the discrete set $\G^a$, we may extend $f$ and $g$ to functions $\widetilde{f}$ and $\widetilde{g}$ in $\ell^2(\G^a)$ by setting them equal to 0 on $\G^x\setminus\G^a_a$. 
This extension provides an isometric embedding of $\ell^2(\G^a_a)$ into $\ell^2(\G^a)$.
The vectors $\widetilde{f}$ and $\widetilde g$ can in turn be extended to continuous sections, say $F$ and $G$, respectively, on the unit space (\cite[Proposition 10.1.10]{Dixmier}). In particular, $F$ and $G$ are continuous sections of the Hilbert bundle associated with the left regular representation of $\G$, as described in \Cref{exa:L2(G)}. 
Moreover, these extensions may be chosen to have the same supremum norm as the original fiber vectors. Indeed, given an extension $F=\{F^b\}_{b\in \G^{(0)}}$ of $\widetilde f$, replace it by the section obtained by multiplying $F^b$ by
$$b\longmapsto \min\left\{1,\frac{\|F^a\|_2}{\|F^b\|_2}\right\},$$
where the multiplier is defined to be 1 whenever $F^b=0$. By \cite[Proposition 10.1.9]{Dixmier}, the resulting vector field remains continuous. It agrees with $F$ at $a$, while its norm is bounded above by $\|F^a\|_2$ at every fiber; hence $\|F\|_\Delta=\sup_{b\in \G^{(0)}}\|F^b\|_2=\|F^a\|_2=\|\widetilde f\|_2$. The same argument applies to $\widetilde g$.

Let $v$ be the coefficient function associated with $F,G$, i.e., $v=\Lambda_{F,G}$. 
For arbitrary $y\in \G^a_a$, we have
$$v(y)=\langle L_y F^a,G^a\rangle_{\ell^2(\G^a)}=\sum_{z\in \G^a}F^a(y^{-1}z)\overline{G^a(z)}=
\sum_{z\in \G^a_a}f(y^{-1}z)\overline{g(z)}=u(y).$$
So $R_a(v)=u$. Since each restriction map is contractive, we get in particular that $\|u\|_{A(G)}\leq \|v\|_{A(\G)}$. On the other hand, $\|v\|_{A(\G)}\leq \|F\|_\Delta\|G\|_\Delta=\|f\|_2\|g\|_2=\|u\|_{A(G)}$. 
So $\|u\|_{A(G)}= \|v\|_{A(\G)}$.
\end{proof}

\begin{prop}\label{prop-spec-etale}
    Let $\G$ be a (locally compact Hausdorff) \'{e}tale groupoid with its counting Haar system. 
    Then $\sigma(A(\G))$ is homeomorphic to $\G$.
\end{prop}

\begin{proof}
Recall that $A(\G)\subseteq C_0(\G)$,
and  $A(\G)$ strongly separates points of $\G$. 
So, the natural embedding $\G\to \sigma(A(\G))$, given by point evaluation, 
is a homeomorphism into $\sigma(A(\G))$ (\cite[Theorem 3.2.4]{Rickart}). It only remains to prove that this embedding is surjective.

Next, consider an arbitrary character $\Phi:A(\G)\to \C$. 
By \Cref{lem-etal-Cc-Acf}, we have $C_c(\G)\subseteq A(\G)$.
We will prove that $\Phi|_{C_c(\G)}$ is contractive, i.e., we prove
\begin{equation}\label{eq-sup-contractive}
    |\Phi(f)|\leq \|f\|_{\sup} \ \text{ for every }\ f\in C_c(G).
\end{equation} 
For $f\in C_c(\G)$ and $c\in \C$ satisfying $|c|>\|f\|_{\sup}$, define
$h:\G\to\C $ as 
$$h(x)=\frac{1}{c}+\frac{f(x)}{c(c-f(x))},$$ 
and note that the summand $\frac{f(x)}{c(c-f(x))}$ belongs to $C_c(\G)$, and thus, to $A(\G)$. 
Let ${\cal A}_1$ denote the unitization of $A(\G)$, and $\widetilde{\Phi}$ be the unique extension of $\Phi$
to a character on ${\cal A}_1$. 
It is easy to see that $h(c1-f)=1$, so $\widetilde{\Phi}(h)\widetilde{\Phi}(c1-f)=1$. 
Therefore, $\widetilde{\Phi}(f)\neq c$, which implies that \eqref{eq-sup-contractive} holds. Consequently, $\Phi$ can be extended to a character on $C_0(\G)$, and it is evaluation at some $x\in \G$.
\end{proof}
\begin{rmk}
    The conclusion of \Cref{prop-spec-etale} holds for any locally compact Hausdorff groupoid $\G$ equipped with a Haar system that satisfies $C_c(\G)\subseteq A(\G)$. In particular, this includes all \'{e}tale groupoids.
\end{rmk}

%%%%%%%%%%%%%%%%%%%%%%%%%%%%%%%%%%%%%%%%%%%%%%%%%%%%%%%%%%%%%%%%%%%%%%%%%%%%%%%%%%%%%%%%%%%%%%%%%%%%%%%%%%%%%%%%%%%%%%%%%%%%%%%%%%%%%%%%%%%%%%

\section*{Acknowledgments}
The first author thank the organizers of the \emph{Northeastern Analysis Meeting (NEAM 2025)}, held at the University of Delaware on November 1--2, 2025, for financial support that enabled JD's visit to the University of Delaware, during which the authors initiated discussions leading to this work.
The second author acknowledges support from National Science Foundation Grant DMS-2408008 and  CCF-2427965 during the preparation of this article. 
The authors acknowledge Simons Foundation Travel Support for Mathematicians that funded their collaboration. 
This material is also based upon work supported by the Swedish Research Council under grant no. 2021-06594 while M.G.~was in residence at Institut Mittag-Leffler in Djursholm, Sweden during the semester ``Operator
algebras and Quantum information theory'' Spring, 2026.
The authors used ChatGPT (OpenAI) to assist with improving the language and presentation of the manuscript in some places. All mathematical content, except \Cref{prop:etale-lift} and \Cref{prop-spec-etale}, was developed and verified by the authors. The statements of \Cref{prop:etale-lift} and \Cref{prop-spec-etale} and an outline of their proof were suggested by ChatGPT. The authors wrote the proofs (in the mathematical style of the paper) and independently verified their correctness.

%\red{We sincerely thank the anonymous reviewers for critically reading the manuscript and suggesting substantial improvements. }
%%%%%%%%%%%%%%%%%%%%%%%%%%%%%%%%%%%%%%%%%%%%%%%%%%%%%%%%%%%%%%%%%%%%%%%%%%%%%%%%%%%%%%%%%%%%%%%%
%%%%%%%%%%%%%%%%%%%%%%%%%%%%%%%%%%%%%%%%%%%%%%%%%%%%%%%%%%%%%%%%%%%%%%%%%%%%%%%%%%%%%%%%%%%%%%%%

\vspace{2ex}
\bibliographystyle{amsplain}
\bibliography{ref}

@article{MuhlyRenaultWilliamsIII,
  author    = {Paul S. Muhly and Jean Renault and Dana P. Williams},
  title     = {Continuous-Trace Groupoid {$C^*$}-Algebras. {III}},
  journal   = {Transactions of the American Mathematical Society},
  volume    = {348},
  number    = {9},
  pages     = {3621--3641},
  year      = {1996},
  doi       = {10.1090/S0002-9947-96-01771-3}
}

@article{GoehleGroupBundle,
  author    = {Geoff Goehle},
  title     = {Group Bundle Duality},
  journal   = {Illinois Journal of Mathematics},
  volume    = {52},
  number    = {3},
  pages     = {951--956},
  year      = {2008},
  doi       = {10.1215/ijm/1254403723}
}

@article{RenaultDualHaar,
  author    = {Renault, Jean},
  title     = {Continuity of the Dual Haar Measure},
  journal   = {Comptes Rendus. Mathématique},
  volume    = {359},
  number    = {4},
  pages     = {415--419},
  year      = {2021},
  doi       = {10.5802/crmath.183}
}

@article{ShirinkalamMedghalchiAmini,
  author    = {Ahmad Shirinkalam and Alireza R. Medghalchi and Massoud Amini},
  title     = {Harmonic Analysis on Locally Compact Abelian Groupoids},
  journal   = {Analele \c{S}tiin\c{t}ifice ale Universit\u{a}\c{t}ii
               ``Alexandru Ioan Cuza'' din Ia\c{s}i. Matematic\u{a}
               (N.S.)},
  volume    = {62},
  number    = {2},
  pages     = {563--570},
  year      = {2016}
}

@misc{EdekoJamneshanKreidler,
  author       = {Nana Edeko and Asgar Jamneshan and Henrik Kreidler},
  title        = {A Peter--Weyl Theorem for Compact Group Bundles and the Geometric Representation of Relatively Ergodic Compact Extensions},
  year         = {2023},
  eprint       = {2302.13630},
  archivePrefix= {arXiv},
  primaryClass = {math.DS}
}

@article {PatersonPaper,
    AUTHOR = {Paterson, Alan L. T.},
     TITLE = {The {F}ourier algebra for locally compact groupoids},
   JOURNAL = {Canad. J. Math.},
  FJOURNAL = {Canadian Journal of Mathematics. Journal Canadien de
              Math\'{e}matiques},
    VOLUME = {56},
      YEAR = {2004},
    NUMBER = {6},
     PAGES = {1259--1289},
      ISSN = {0008-414X},
   MRCLASS = {43A32 (46J99)},
  MRNUMBER = {2102633},
MRREVIEWER = {Humberto E. Prado},
       DOI = {10.4153/CJM-2004-055-8},
       URL = {https://doi-org.udel.idm.oclc.org/10.4153/CJM-2004-055-8},
}

@book {Dixmier,
    AUTHOR = {Dixmier, Jacques},
     TITLE = {{$C\sp*$}-algebras},
    SERIES = {North-Holland Mathematical Library, Vol. 15},
      NOTE = {Translated from the French by Francis Jellett},
 PUBLISHER = {North-Holland Publishing Co., Amsterdam-New York-Oxford},
      YEAR = {1977},
     PAGES = {xiii+492},
      ISBN = {0-7204-0762-1},
   MRCLASS = {46L05},
  MRNUMBER = {0458185},
}

@article {Eymard,
    AUTHOR = {Eymard, Pierre},
     TITLE = {L'alg\`ebre de {F}ourier d'un groupe localement compact},
   JOURNAL = {Bull. Soc. Math. France},
  FJOURNAL = {Bulletin de la Soci\'{e}t\'{e} Math\'{e}matique de France},
    VOLUME = {92},
      YEAR = {1964},
     PAGES = {181--236},
      ISSN = {0037-9484},
   MRCLASS = {22.65},
  MRNUMBER = {228628},
MRREVIEWER = {C. Dunkl},
       URL = {http://www.numdam.org/item?id=BSMF_1964__92__181_0},
}

@article {Renault,
    AUTHOR = {Renault, Jean},
     TITLE = {The {F}ourier algebra of a measured groupoid and its
              multipliers},
   JOURNAL = {J. Funct. Anal.},
  FJOURNAL = {Journal of Functional Analysis},
    VOLUME = {145},
      YEAR = {1997},
    NUMBER = {2},
     PAGES = {455--490},
      ISSN = {0022-1236},
   MRCLASS = {43A15 (22A22 43A22 43A35 46L99)},
  MRNUMBER = {1444088},
MRREVIEWER = {Fran\c{c}oise Lust-Piquard},
       DOI = {10.1006/jfan.1996.3039},
       URL = {https://doi.org/10.1006/jfan.1996.3039},
}

@book {Rickart,
    AUTHOR = {Rickart, Charles E.},
     TITLE = {General theory of {B}anach algebras},
    SERIES = {The University Series in Higher Mathematics},
 PUBLISHER = {D. Van Nostrand Co., Inc., Princeton, N.J.-Toronto-London-New
              York},
      YEAR = {1960},
     PAGES = {xi+394},
   MRCLASS = {46.00},
  MRNUMBER = {0115101},
MRREVIEWER = {M. Nakamura},
}

@article {Willett,
    AUTHOR = {Willett, Rufus},
     TITLE = {A non-amenable groupoid whose maximal and reduced
              {$C^*$}-algebras are the same},
   JOURNAL = {M\"{u}nster J. Math.},
  FJOURNAL = {M\"{u}nster Journal of Mathematics},
    VOLUME = {8},
      YEAR = {2015},
    NUMBER = {1},
     PAGES = {241--252},
      ISSN = {1867-5778},
   MRCLASS = {46L55 (22A22 22D25 46L05)},
  MRNUMBER = {3549528},
MRREVIEWER = {Jean N. Renault},
       DOI = {10.17879/65219671638},
       URL = {https://doi.org/10.17879/65219671638},
}

@article {Skandalis,
    AUTHOR = {Higson, Nigel and Lafforgue, Vincent and Skandalis, Georges},
     TITLE = {Counterexamples to the {B}aum-{C}onnes conjecture},
   JOURNAL = {Geom. Funct. Anal.},
  FJOURNAL = {Geometric and Functional Analysis},
    VOLUME = {12},
      YEAR = {2002},
    NUMBER = {2},
     PAGES = {330--354},
      ISSN = {1016-443X},
   MRCLASS = {19K56 (22A22 46L80 46L85 58J22)},
  MRNUMBER = {1911663},
MRREVIEWER = {Alain Valette},
       DOI = {10.1007/s00039-002-8249-5},
       URL = {https://doi.org/10.1007/s00039-002-8249-5},
}

@article {Leptin,
    AUTHOR = {Leptin, Horst},
     TITLE = {Sur l'alg\`ebre de {F}ourier d'un groupe localement compact},
   JOURNAL = {C. R. Acad. Sci. Paris S\'{e}r. A-B},
  FJOURNAL = {Comptes Rendus Hebdomadaires des S\'{e}ances de l'Acad\'{e}mie des
              Sciences. S\'{e}ries A et B},
    VOLUME = {266},
      YEAR = {1968},
     PAGES = {A1180--A1182},
      ISSN = {0151-0509},
   MRCLASS = {22.65 (46.00)},
  MRNUMBER = {239002},
MRREVIEWER = {M. Rajagopalan},
}

@article {Brandt,
    AUTHOR = {Brandt, Heinrich},
     TITLE = {\"{U}ber eine {V}erallgemeinerung des {G}ruppenbegriffes},
   JOURNAL = {Math. Ann.},
  FJOURNAL = {Mathematische Annalen},
    VOLUME = {96},
      YEAR = {1927},
    NUMBER = {1},
     PAGES = {360--366},
      ISSN = {0025-5831,1432-1807},
   MRCLASS = {99-04},
  MRNUMBER = {1512323},
       DOI = {10.1007/BF01209171},
       URL = {https://doi.org/10.1007/BF01209171},
}

@book {PatersonBook,
    AUTHOR = {Paterson, Alan L. T.},
     TITLE = {Groupoids, inverse semigroups, and their operator algebras},
    SERIES = {Progress in Mathematics},
    VOLUME = {170},
 PUBLISHER = {Birkh\"{a}user Boston, Inc., Boston, MA},
      YEAR = {1999},
     PAGES = {xvi+274},
      ISBN = {0-8176-4051-7},
   MRCLASS = {22A22 (20M18 22-02 22D25 46L05)},
  MRNUMBER = {1724106},
MRREVIEWER = {Jean N. Renault},
       DOI = {10.1007/978-1-4612-1774-9},
       URL = {https://doi.org/10.1007/978-1-4612-1774-9},
}

@article {Brown-etal-2014,
    AUTHOR = {Brown, Jonathan and Clark, Lisa Orloff and Farthing, Cynthia
              and Sims, Aidan},
     TITLE = {Simplicity of algebras associated to \'etale groupoids},
   JOURNAL = {Semigroup Forum},
  FJOURNAL = {Semigroup Forum},
    VOLUME = {88},
      YEAR = {2014},
    NUMBER = {2},
     PAGES = {433--452},
      ISSN = {0037-1912,1432-2137},
   MRCLASS = {46L05 (22A22)},
  MRNUMBER = {3189105},
MRREVIEWER = {Jean\ N.\ Renault},
       DOI = {10.1007/s00233-013-9546-z},
       URL = {https://doi.org/10.1007/s00233-013-9546-z},
}

@article {Pedro,
    AUTHOR = {Resende, Pedro},
     TITLE = {\'Etale groupoids and their quantales},
   JOURNAL = {Adv. Math.},
  FJOURNAL = {Advances in Mathematics},
    VOLUME = {208},
      YEAR = {2007},
    NUMBER = {1},
     PAGES = {147--209},
      ISSN = {0001-8708,1090-2082},
   MRCLASS = {22A22 (18B40 54B30 54H10)},
  MRNUMBER = {2304314},
MRREVIEWER = {R.\ Gylys},
       DOI = {10.1016/j.aim.2006.02.004},
       URL = {https://doi.org/10.1016/j.aim.2006.02.004},
}

@book {ToolKit,
    AUTHOR = {Williams, Dana P.},
     TITLE = {A tool kit for groupoid {$C^*$}-algebras},
    SERIES = {Mathematical Surveys and Monographs},
    VOLUME = {241},
 PUBLISHER = {American Mathematical Society, Providence, RI},
      YEAR = {2019},
     PAGES = {xv+398},
      ISBN = {978-1-4704-5133-2},
   MRCLASS = {46-02 (22A22 22D25 46Lxx)},
  MRNUMBER = {3969970},
MRREVIEWER = {Valentin\ Deaconu},
       DOI = {10.1016/j.physletb.2019.06.021},
       URL = {https://doi.org/10.1016/j.physletb.2019.06.021},
}

@article {Oty,
    AUTHOR = {Oty, Karla J.},
     TITLE = {Fourier-{S}tieltjes algebras of {$r$}-discrete groupoids},
   JOURNAL = {J. Operator Theory},
  FJOURNAL = {Journal of Operator Theory},
    VOLUME = {41},
      YEAR = {1999},
    NUMBER = {1},
     PAGES = {175--197},
      ISSN = {0379-4024},
   MRCLASS = {43A32},
  MRNUMBER = {1675192},
}

@article {Herz,
    AUTHOR = {Herz, Carl},
     TITLE = {Harmonic synthesis for subgroups},
   JOURNAL = {Ann. Inst. Fourier (Grenoble)},
  FJOURNAL = {Universit\'{e} de Grenoble. Annales de l'Institut Fourier},
    VOLUME = {23},
      YEAR = {1973},
    NUMBER = {3},
     PAGES = {91--123},
      ISSN = {0373-0956,1777-5310},
   MRCLASS = {43A45},
  MRNUMBER = {355482},
MRREVIEWER = {D.\ E.\ Ramirez},
       URL = {http://www.numdam.org/item?id=AIF_1973__23_3_91_0},
}

@article {Cowling,
    AUTHOR = {Cowling, Michael and Rodway, Paul},
     TITLE = {Restrictions of certain function spaces to closed subgroups of
              locally compact groups},
   JOURNAL = {Pacific J. Math.},
  FJOURNAL = {Pacific Journal of Mathematics},
    VOLUME = {80},
      YEAR = {1979},
    NUMBER = {1},
     PAGES = {91--104},
      ISSN = {0030-8730,1945-5844},
   MRCLASS = {43A30 (43A60)},
  MRNUMBER = {534697},
MRREVIEWER = {C.\ Chou},
       URL = {http://projecteuclid.org/euclid.pjm/1102785953},
}

@article {MahyaThesisPaper,
    AUTHOR = {Ghandehari, Mahya},
     TITLE = {Amenability properties of {R}ajchman algebras},
   JOURNAL = {Indiana Univ. Math. J.},
  FJOURNAL = {Indiana University Mathematics Journal},
    VOLUME = {61},
      YEAR = {2012},
    NUMBER = {3},
     PAGES = {1369--1392},
      ISSN = {0022-2518,1943-5258},
   MRCLASS = {43A30 (43A07 46J10 47B47)},
  MRNUMBER = {3071703},
MRREVIEWER = {Maedeh\ Soroushmehr},
       DOI = {10.1512/iumj.2012.61.4622},
       URL = {https://doi.org/10.1512/iumj.2012.61.4622},
}

@book {HewittII,
    AUTHOR = {Hewitt, Edwin and Ross, Kenneth A.},
     TITLE = {Abstract harmonic analysis. {V}ol. {II}: {S}tructure and
              analysis for compact groups. {A}nalysis on locally compact
              {A}belian groups},
    SERIES = {Die Grundlehren der mathematischen Wissenschaften},
    VOLUME = {Band 152},
 PUBLISHER = {Springer-Verlag, New York-Berlin},
      YEAR = {1970},
     PAGES = {ix+771},
   MRCLASS = {42.50 (22.00)},
  MRNUMBER = {262773},
MRREVIEWER = {R.\ E.\ Edwards},
}

@article{DeG,
    author = {DeGaetani, Joseph and Ghandehari, Mahya},
    title = {On the Restriction Maps of the Fourier and Fourier–Stieltjes Algebras Over Locally Compact Groupoids},
    journal = {International Mathematics Research Notices},
    volume = {2025},
    number = {6},
    pages = {rnaf057},
    year = {2025},
    month = {03},
    issn = {1073-7928},
    doi = {10.1093/imrn/rnaf057},
    url = {https://doi.org/10.1093/imrn/rnaf057},
    eprint = {https://academic.oup.com/imrn/article-pdf/2025/6/rnaf057/62427575/rnaf057.pdf},
}

@article {MR4645718,
    AUTHOR = {Jiang, Baojie and Zhang, Jiawen and Zhang, Jianguo},
     TITLE = {Quasi-locality for \'{e}tale groupoids},
   JOURNAL = {Comm. Math. Phys.},
  FJOURNAL = {Communications in Mathematical Physics},
    VOLUME = {403},
      YEAR = {2023},
    NUMBER = {1},
     PAGES = {329--379},
      ISSN = {0010-3616,1432-0916},
   MRCLASS = {46L05 (46L80)},
  MRNUMBER = {4645718},
       DOI = {10.1007/s00220-023-04782-x},
       URL = {https://doi.org/10.1007/s00220-023-04782-x},
}

@article {MR4608439,
    AUTHOR = {Kranz, Julian},
     TITLE = {Amenability for actions of \'{e}tale groupoids on
              {$C^*$}-algebras and {F}ell bundles},
   JOURNAL = {Trans. Amer. Math. Soc.},
  FJOURNAL = {Transactions of the American Mathematical Society},
    VOLUME = {376},
      YEAR = {2023},
    NUMBER = {7},
     PAGES = {5089--5121},
      ISSN = {0002-9947,1088-6850},
   MRCLASS = {46L55 (22A22)},
  MRNUMBER = {4608439},
MRREVIEWER = {Eusebio\ Gardella},
       DOI = {10.1090/tran/8909},
       URL = {https://doi.org/10.1090/tran/8909},
}

@article {MR4592883,
    AUTHOR = {Christensen, Johannes},
     TITLE = {The structure of {KMS} weights on \'{e}tale groupoid
              {$C^*$}-algebras},
   JOURNAL = {J. Noncommut. Geom.},
  FJOURNAL = {Journal of Noncommutative Geometry},
    VOLUME = {17},
      YEAR = {2023},
    NUMBER = {2},
     PAGES = {663--691},
      ISSN = {1661-6952,1661-6960},
   MRCLASS = {46L05 (37A55 46L30)},
  MRNUMBER = {4592883},
       DOI = {10.4171/jncg/507},
       URL = {https://doi.org/10.4171/jncg/507},
}

@article {MR4563262,
    AUTHOR = {Neshveyev, Sergey and Schwartz, Gaute},
     TITLE = {Non-{H}ausdorff \'{e}tale groupoids and {${\rm C}^*$}-algebras
              of left cancellative monoids},
   JOURNAL = {M\"{u}nster J. Math.},
  FJOURNAL = {M\"{u}nster Journal of Mathematics},
    VOLUME = {16},
      YEAR = {2023},
    NUMBER = {1},
     PAGES = {147--175},
      ISSN = {1867-5778,1867-5786},
   MRCLASS = {46L05 (20M05 46L55)},
  MRNUMBER = {4563262},
MRREVIEWER = {Qingzhai\ Fan},
}

@article {MR4533452,
    AUTHOR = {Bice, Tristan},
     TITLE = {Representing structured semigroups on \'{e}tale groupoid
              bundles},
   JOURNAL = {J. Aust. Math. Soc.},
  FJOURNAL = {Journal of the Australian Mathematical Society},
    VOLUME = {114},
      YEAR = {2023},
    NUMBER = {1},
     PAGES = {1--49},
      ISSN = {1446-7887,1446-8107},
   MRCLASS = {20M18 (06F05 20M30 22A22 46L05 54D80)},
  MRNUMBER = {4533452},
       DOI = {10.1017/s1446788722000064},
       URL = {https://doi.org/10.1017/s1446788722000064},
}

@book {MR4510931,
    AUTHOR = {Exel, Ruy and Pitts, David R.},
     TITLE = {Characterizing groupoid {$\rm C^*$}-algebras of
              non-{H}ausdorff \'{e}tale groupoids},
    SERIES = {Lecture Notes in Mathematics},
    VOLUME = {2306},
 PUBLISHER = {Springer, Cham},
      YEAR = {[2022] \copyright 2022},
     PAGES = {viii+156},
      ISBN = {978-3-031-05512-6; 978-3-031-05513-3},
   MRCLASS = {46-02 (22A22 46L05 46L45 46L55)},
  MRNUMBER = {4510931},
MRREVIEWER = {Cristian\ Ivanescu},
       DOI = {10.1007/978-3-031-05513-3},
       URL = {https://doi.org/10.1007/978-3-031-05513-3},
}

@article {MR4483992,
    AUTHOR = {Matui, Hiroki},
     TITLE = {Long exact sequences of homology groups of \'{e}tale
              groupoids},
   JOURNAL = {Discrete Contin. Dyn. Syst.},
  FJOURNAL = {Discrete and Continuous Dynamical Systems. Series A},
    VOLUME = {42},
      YEAR = {2022},
    NUMBER = {11},
     PAGES = {5239--5271},
      ISSN = {1078-0947,1553-5231},
   MRCLASS = {55U35 (19L99 37B52)},
  MRNUMBER = {4483992},
MRREVIEWER = {Jes\'{u}s\ Rodr\'{\i}guez-L\'{o}pez},
       DOI = {10.3934/dcds.2022095},
       URL = {https://doi.org/10.3934/dcds.2022095},
}

@article {Matui-PAMS-2012,
    AUTHOR = {Matui, Hiroki},
     TITLE = {Homology and topological full groups of \'etale groupoids on
              totally disconnected spaces},
   JOURNAL = {Proc. Lond. Math. Soc. (3)},
  FJOURNAL = {Proceedings of the London Mathematical Society. Third Series},
    VOLUME = {104},
      YEAR = {2012},
    NUMBER = {1},
     PAGES = {27--56},
      ISSN = {0024-6115,1460-244X},
   MRCLASS = {19D55 (22A22 46L55)},
  MRNUMBER = {2876963},
MRREVIEWER = {Hiroyuki\ Osaka},
       DOI = {10.1112/plms/pdr029},
       URL = {https://doi.org/10.1112/plms/pdr029},
}

@book{dissertation,
  title={Banach Algebraic Properties of the Fourier and Fourier-Stieltjes Algebras over Locally Compact Groupoids},
  author={DeGaetani, Joseph},
  year={2024},
  publisher={University of Delaware}
}

@article {Ramsay-Martin-97, 
    AUTHOR = {Ramsay, Arlan and Walter, Martin E.},
     TITLE = {Fourier-{S}tieltjes algebras of locally compact groupoids},
   JOURNAL = {J. Funct. Anal.},
  FJOURNAL = {Journal of Functional Analysis},
    VOLUME = {148},
      YEAR = {1997},
    NUMBER = {2},
     PAGES = {314--367},
      ISSN = {0022-1236,1096-0783},
   MRCLASS = {22D25 (43A40 46L05)},
  MRNUMBER = {1469344},
MRREVIEWER = {Jean\ N.\ Renault},
       DOI = {10.1006/jfan.1996.3083},
       URL = {https://doi.org/10.1006/jfan.1996.3083},
}

@incollection {Paterson-2003,
    AUTHOR = {Paterson, Alan L. T.},
     TITLE = {The {F}ourier-{S}tieltjes and {F}ourier algebras for locally
              compact groupoids},
 BOOKTITLE = {Trends in {B}anach spaces and operator theory ({M}emphis,
              {TN}, 2001)},
    SERIES = {Contemp. Math.},
    VOLUME = {321},
     PAGES = {223--237},
 PUBLISHER = {Amer. Math. Soc., Providence, RI},
      YEAR = {2003},
      ISBN = {0-8218-3234-4},
   MRCLASS = {43A20 (22A22 43A25 43A35 46L08 46L87)},
  MRNUMBER = {1978820},
       DOI = {10.1090/conm/321/05647},
       URL = {https://doi.org/10.1090/conm/321/05647},
}
\vspace{-0.75ex}

\end{document}